\documentclass[11pt,a4paper]{article}
\usepackage{amsmath}
\usepackage{amsfonts}
\usepackage{amssymb}
\usepackage{graphicx}
\usepackage[utf8]{inputenc}
\usepackage{array}
\usepackage{t1enc}
\usepackage[mathscr]{eucal}
\usepackage{amsthm}
\usepackage{mathtools}

\usepackage[usenames,dvipsnames]{color}

\usepackage[colorlinks=false,pdfborderstyle={/W 1}]{hyperref}
\def\arXiv#1#2{\href{http://front.math.ucdavis.edu/#1}{{\tt arXiv:#1 [#2]}}} 
\def\arXivo#1{\href{http://front.math.ucdavis.edu/#1}{{\tt [arXiv:#1]}}}

\numberwithin{equation}{section}
\numberwithin{table}{section}
\numberwithin{figure}{section}
\newenvironment*{pf}{\noindent{\it Proof.}}{\hfill $\blacksquare$}
\newenvironment*{pfo}{\medskip\noindent{\it Proof of}}{\hfill $\blacksquare$}

  \ifx\LabelFigloaded\MYundefined\relax
  \else
   \fi

  \def\LabelFigloaded{\relax}% now loaded

  \chardef\LabelFigCatAt\the\catcode`\@
  \catcode`\@=11

 \let\LabelFigwlog@ld\wlog
 \def\wlog#1{\relax}

 \ifx\\\MYundefined@
    \let\\\relax
 \fi

  \def\ms@g{\immediate\write16}

 \def\N@wif{\csname newif\endcsname }
 \def\Temp@ {\N@wif\ifIN@}
 \ifx\INN@\MYundefined@
    \else \let\Temp@\relax
 \fi
 \Temp@

  \def\IN@{\expandafter\INN@\expandafter}
  \long\def\INN@0#1@#2@{\long\def\NI@##1#1##2##3\ENDNI@
    {\ifx\m@rker##2\IN@false\else\IN@true\fi}%
     \expandafter\NI@#2@@#1\m@rker\ENDNI@}
  \def\m@rker{\m@@rker}
 
  \newtoks\Initialtoks@  \newtoks\Terminaltoks@
  \def\SPLIT@{\expandafter\SPLITT@\expandafter}
  \def\SPLITT@0#1@#2@{\def\TTILPS@##1#1##2@{%
     \Initialtoks@{##1}\Terminaltoks@{##2}}\expandafter\TTILPS@#2@}

 \def\Shifted@@#1#2#3{\setbox0=\hbox{#3}%
   \raise -\dp0\vbox {\kern-#2%
       \hbox {\kern#1\unhbox0\kern-#1}%
           \kern#2}}

 \newcount\gridcount
 \newbox\auxGridbox@ \newbox\hGridbox@ \newbox\vGridbox@
 \newbox\Labelbox@ \newbox\auxLabelbox@
 \newbox\Coordinatebox@
 \newtoks\Labeltoks@
 \newdimen\Wdd@ \newdimen\Htt@
 \newdimen\Wddd@ \newdimen\Httt@
 
 \def\Wr@{\immediate\write16}

 \newdimen\GL@wd%% grid-line width
 \def\GridLineWidth#1{\GL@wd=#1}

 \def\gobble#1{}
 \def\EdgeErr@{\Wr@{}%
      \Wr@{\string\Edges\space argument
      1, 10, 100 or 1000 please\string!}%
      }

 \newcount\Edgect@

 \def\Sweepup#1\endSweepup{}

 \def\SetEdges@{%
    \edef\Zr@@s{\expandafter\gobble\number\Edgect@\empty}%
        \count255=0\Zr@@s\relax
        \ifnum\count255=\z@\else\EdgeErr@\show\tailtest\fi
        \count255=1\Zr@@s\relax%\showthe\count255
        \ifnum\count255=\Edgect@\relax\else\EdgeErr@\show\leadtest\fi
    \EdgGl@b\edef\Zr@s{\expandafter\gobble\Zr@@s\empty}%\show\Zr@s
    \ifnum\Edgect@>\@ne\relax\EdgGl@b\let\L@Dc\empty
        \else\EdgGl@b\edef\L@Dc{\string.}\fi
    \ifnum\Edgect@>\@ne\relax
        \EdgGl@b\edef\Edgescale@##1{\divide##1 by \Edgect@}%
        \else\EdgGl@b\edef\Edgescale@##1{}\fi
    }

 \def\Edges#1{\Edgect@=#1\relax
     \let\EdgGl@b\global \SetEdges@}

 \Edges{1}%% default

 \def\hhrule{\hrule height \GL@wd\vskip-.\GL@wd}

 \def\hRule@{%
   \advance\gridcount -2%
   \vfil\hhrule\vfil
   \llap{\smash{\raise -2.5pt
     \hbox{\L@Dc\number\gridcount\Zr@s\kern2pt}}}%
   \hhrule
   }

\def\vvrule{\vrule width \GL@wd \kern-\GL@wd}

 \def\vRule@{\advance\gridcount 2%
   \hfil\vvrule\hfil
   \setbox\auxGridbox@=\vbox to 0pt
      {\vskip \Htt@\vskip 2pt
        \hbox to 0pt{\hss\L@Dc\number\gridcount\Zr@s\hss}\vss}%
      \wd\auxGridbox@=0pt \box\auxGridbox@
   \vvrule
   }

 \def\PlaceGrid@@{\gridcount=10 
  \setbox\hGridbox@=\hbox{%
        \hbox{%
             \hskip-.4pt\vrule
             \vbox to \Htt@{%
               \offinterlineskip\parindent=\z@\relax
               \hbox to \Wdd@{\hfil}
               \hRule@\hRule@\hRule@\hRule@
               \vfil\hhrule\vfil}%
             \vrule\hskip-.4pt}
    }%
  \gridcount=0%
  \setbox\vGridbox@=\hbox{%
      \vbox{\offinterlineskip\parindent=0pt\hsize=0pt
         \vskip-.4pt\hrule%
         \hbox to \Wdd@{%
                 \vtop to \Htt@{\vfil}%
                 \vRule@\vRule@\vRule@\vRule@
                 \hfil\vvrule\hfil}%
         \hrule\vskip-.4pt}}%
  \wd\hGridbox@=0pt\ht\hGridbox@=0pt
  \wd\vGridbox@=0pt\ht\vGridbox@=0pt
  \hbox{\box\hGridbox@\box\vGridbox@}%
  }

 \def\LabelsGlobal{\def\LabGl@b{\global}}
 \def\LabelsLocal{\def\LabGl@b{}}
 \LabelsGlobal %% default

 \def\SetLabels#1\endSetLabels{%
   \LabGl@b\Labeltoks@={#1()\\}%
   }

 \LabGl@b\Labeltoks@={()\\}

 \def\ShowGrid{\LabGl@b\let\PlaceGrid@\PlaceGrid@@}
 \def\HideGrid{\LabGl@b\let\PlaceGrid@\relax}
 \def\Grids{\ShowGrid\LabGl@b\let\GridSwitch@\ShowGrid}
 \def\noGrids{\HideGrid\LabGl@b\let\GridSwitch@\HideGrid}

 \noGrids

 \def\bAdjust@@{%
     \setbox\auxLabelbox@=\hbox{\raise \dp\auxLabelbox@
            \box\auxLabelbox@}}
 \def\bAdjust@{\let\vAdjust@\bAdjust@@}

 \def\eAdjust@@{\dimen0=-.5\ht\auxLabelbox@
     \advance\dimen0 by .5\dp\auxLabelbox@
     \setbox\auxLabelbox@=
            \hbox{\raise\dimen0\box\auxLabelbox@}}
 \def\eAdjust@{\let\vAdjust@\eAdjust@@}

 \def\tAdjust@@{%
     \setbox\auxLabelbox@=\hbox{\raise-\ht\auxLabelbox@
            \box\auxLabelbox@}}
 \def\tAdjust@{\let\vAdjust@\tAdjust@@}

 \let\vAdjust@\relax

 \def\lAdjust@{\let\hAdjust@\rlap}
 \def\rAdjust@{\let\hAdjust@\llap}

 \let\hAdjust@\relax\let\vAdjust@\relax

 \def\FetchLabel@#1(#2)#3\\{%
     \IN@0#2@@\ifIN@
        \setbox0=\hbox{\ignorespaces#1#3\unskip}%
        \ifdim\wd0>0pt
           \ms@g{}%
           \ms@g{ !!! Bad label(s)? !!!}%
           \message{ #1(#2)#3}%
        \fi
        \def\LabelMole@##1\endFetchLabel@{%
            \IN@0()\\@##1@%
            \ifIN@\def\Temp@{\FetchLabel@##1\endFetchLabel@}%
            \else\def\Temp@{}%
            \fi
            \Temp@
           }%
     \else
       \ignorespaces#1\unskip
       \setbox\auxLabelbox@=%
         \hbox to 0pt{\hss\ignorespaces\hAdjust@
          {\ignorespaces#3\unskip}\hss}%
       \vAdjust@
       \let\hAdjust@\relax\let\vAdjust@\relax
       \AugmentLabelBox@@{#2}%
       \ht\Labelbox@=0pt\dp\Labelbox@=0pt
       \let\LabelMole@\FetchLabel@%
     \fi\LabelMole@}

 \newtoks\XYSep@ %\XYSep@{*}
 \def\SetXYSeparator#1{%
     \IN@0#1@@\ifIN@\XYSep@{*}%
     \else
     \XYSep@{#1}%
     \fi
     }

 \SetXYSeparator*

 \def\AugmentLabelBox@@#1{%
     \IN@0\the\XYSep@ @#1@\ifIN@
       \SPLIT@0\the\XYSep@ @#1@%
       \setbox\Labelbox@=\hbox to 0pt{%
         \unhbox\Labelbox@
         \Shifted@@{\the\Initialtoks@\Wddd@}%
         {\the\Terminaltoks@\Httt@}%
         {\box\auxLabelbox@}}%
     \else
         \ms@g{}%
         \ms@g{ !!! Bad insertion point. !!!}%
         \message{ (#1\ this point was rejected.)}%
     \fi
    }

 \def\FetchOption@#1[#2]#3\endFetchOption@{%
    \def\temp{#1}%\show\temp
    \ifx\temp\empty
       \Edgect@=#2\relax%\showthe\Edgect@
       \let\EdgGl@b\relax
       \SetEdges@%\def\Edgescale@##1{\divide##1 by \Edgect@\relax}%
       \Cleaner@#3%
    \fi}

 \def\Cleaner@#1[@]{\Labeltoks@{#1}}
     
 \def\PlaceLabels@@{\mathsurround=0pt%\bgroup
     \def\Cr@{\\}%
     \let\L\lAdjust@\let\R\rAdjust@
     \let\B\bAdjust@\let\E\eAdjust@\let\T\tAdjust@
     \expandafter\FetchOption@\the\Labeltoks@[@]\endFetchOption@
     \Wddd@=\Wdd@ \Edgescale@\Wddd@ %\showthe\Edgect@
     \Httt@=\Htt@ \Edgescale@\Httt@
     \expandafter\FetchLabel@\the\Labeltoks@\endFetchLabel@
     \box\Labelbox@%\egroup
     }%

 \let \PlaceLabels@\PlaceLabels@@

 \def\AffixLabels#1{\setbox\Coordinatebox@=\hbox{#1}%
      \Wdd@=\wd\Coordinatebox@ \Htt@=\ht\Coordinatebox@
      \advance\Htt@ \dp\Coordinatebox@
      \hbox{\copy\Coordinatebox@\kern-\Wdd@ 
           \Shifted@@{0pt}{-\dp\Coordinatebox@}%
           {\PlaceLabels@\PlaceGrid@}%
           \kern\Wdd@}%
      \GridSwitch@ %% next grid hidden
      \LabGl@b\Labeltoks@{()\\}%
      }
 
   \let\wlog\LabelFigwlog@ld   %%restore logging
   \catcode`\@=\LabelFigCatAt  %%12 or 13

                                By

              Raymond S\'eroul <A18645@FRCCSC21.BITNET>
                                and 
              Laurent Siebenmann <lcs@topo.math.u-psud.fr>
    
              VERSIONS: July 1991, Oct 1991, Jan 1992, July 1992

INTRODUCTION

      This labelling package is intended for TeX users who
rely on non-TeX sources for for their graphics inserts.  It
provides means for adding TeX labels to such inserts with a
minimum of fuss. 

       For most labels, TeX users have in the past found it
reasonably convenient to rely on non-TeX sources. Typical
occasions when an inescapable need for TeX labels seemed to
arise are

 (a) when the graphics program lacks certain exotic or complex
mathematical symbols

 (b) when the very highest typographical quality is wanted for the
labels

 (c) when labels included with the graphics fail to print, 
 and you cannot figure out why (cf. boxedeps.doc).  The labels
 provided by labelfig.tex are 100% portable.

       Since this package first appeared, many users, who in the
past scarcely dreamed of using TeX labels, have come to use
nothing but.  So it is now appropriate to add

Intoxication Warning:  TeX labels may be addictive and expensive. 

     If you have a fast preview you may disagree, and even find
that this package provides an agreeable paste-up environment; see
extra applications at end.

     Note to publishers: It is possible and convenient to ultimately
export the TeX labels produced by labelfig.tex to become an integral
part of the EPS file. This is often desired by a publisher who typically
uses an "upmarket" graphics or page layout program, with which the
staff is skilled in perfecting figures.  See Appendix I for
a recipe.

     The authors are grateful to Patrick Ion of Math Reviews for
helpful comments and encouragement.

BASIC INSTRUCTIONS

    After reading in the macro file using

preview or proof your figure with a coordinate grid printed on
top, by typing the following:

    \ShowGrid  % shows grid  for next figure only
    \AffixLabels{<the graphics insertion>}

Here <the graphics insertion> is what you would type to insert
the graphics object alone without the grid.  This must provide
for the space around it. For example <the graphics insertion>
might well be \BoxedEPSF{MyFigure scaled 700} using the
boxedeps.tex macro package (from same source); this provides a
TeX box containing the encapsulated PostScript insert specified by
the file MyFigure. \AffixLabels{...} provides the grid (supposing
\ShowGrid is present) and later, once you have specified labels
using the grid, it will "tack on" the labels.

     The grid is a sort of (usually elongated) checkerboard of
ten rows and ten columns and its (internal) partitions are by
default numbered  .1, ... ,.9  both horizontally (X-coordinate
running left to right) and vertically (Y-coordinate running bottom
to top).  Thus the points enclosed by the grid correspond to the
points of the unit square in the cartesian "X-Y" plane, the lower
left corner corresponding to the origin (0,0).  By extrapolation,
the full page corresponds to a larger rectangle in the plane.

     These coordinates serve to position labels as follows.
Before the \AffixLabels{...} command type label specifications:

  \SetLabels
   (<X-coordinate>*<Y-coordinate>) <first label> \\
   .
   .
   .
   (<X-coordinate>*<Y-coordinate>)  <last label> \\
  \endSetLabels

Each row specifies one label and is terminated by \\.  In each
row, the position indicator comes first; it is written as a
standard cartesian point except that the X- and Y- coordinates
are separated by * rather than a comma because TeX allows a
comma as decimal point. There are no dimension units to specify
as the unit is the grid itself.

     By default, this cartesian point specifies where the middle
of the baseline of the label will be located.  However if you precede
the point by \L [or \R] the left [or right] edge of the baseline will
be located there. Similarly you may also precede the point by \T, \E,
or \B to vertically align the top equator or bottom of the label box
at the specified point.  This gives nine standard positions of
the label with respect to the insertion point --- corresponding to
the eight principle points of the compas and the center

                     \L\T     \T      \R\T

                     \L\E     \E      \R\E

                     \L\B     \B      \R\B

But this neglects the default "baseline" level of TeX,
giving potentially three more positions

                     \L    <no tag>   \R

For text, the baseline level is often the preferred. Its relation to
the others is variable. It will often coincide with the bottom level,
as happens for "X".  But it is often distinct, as for "g", in which
case you have in all 12 distinct positions rather than 9.

     It is convenient to think of this specification of label
position as attaching the label by a thumb-tack to the coordinate
grid. There are up to twelve positions of the thumb-tack on the
label, while the position of the thumb-tack on the coordinate grid is
arbitrary.  Normally, one choses the position of the thumb-tack on
the label to be the one that is the closest to the item being
labeled.  There are good reasons for this "rule of thumb":

   (a)  It facilitates correct positioning at first try.

   (b)  If the scale of the figure must be altered after labels
have been affixed, the labels have a good chance of remaining well
positioned.

   (c)  The visible grid need not extend beyond the "bounding box"
for the figure, because the best preferred position is always
(at least almost) within the bounding box .

The second reason is particularly important. Indeed it often
happens that scale has to be altered after labelling begins, in
order to either provide space for the labels, or to adjust
proportions between the labels and the figure.  (The size of labels
is unaffected by scaling.)

     Here is an artificial but self-contained test which uses
TeX rules to make a graphics object.

TEST

    Do not skip this!

 \def\FrameIt#1{\hbox{\vrule$\vcenter {\hrule\kern3pt%
             \hbox {\kern3pt #1\kern3pt}%
               \kern3pt\hrule}$\relax\vrule}}

 \def\Caption#1#2{\FrameIt{%
       \vtop {\hsize=#1\relax \parindent=0pt
         \leftskip=0pt \rightskip=0pt plus15pt
         \parfillskip=0pt
         \lineskip=1pt\baselineskip=0pt
         #2}}}

 \def\FirstQuadrant{\hbox to 100pt{\vrule\vbox to 100pt{%
        \hbox to 100pt{\hfil}\vfil\hrule}\hss}}

  \SetLabels
    \R(.5*.2) $\zeta\,\cdot$\\
    (.9*-.10) $\xi$\\
    \R(-.03*.9) $\eta$\\
    \T(.5*.9) \Caption{70pt}{%
          \it The norm of
          $g(\xi+i\eta)$ is indicated on
          contours of this invisible surface.}\\
  \endSetLabels

  \AffixLabels{\FirstQuadrant} \end

  Note that the coordinates to use for labels are indicated on the
edges of the grid (when visible) corresponding to the conventional
x- and y- axes of the Cartesian plane. By default the grid is
1-by-1. However, by the command \Edges{100}, you can change this
to 100-by-100 and many users find this alternative most
convenient. Place the command \Edges{...} in your style file (or
header) since its effect is is global. Other possible edge values
are 10 and 1000.

  If you use the command \Edges{...} at all, do so with care.  For
if you accidentally delete an \Edges{...} command your labels will
abruptly be badly misplaced and may logically but mysteriously
generate "dimension too big" errors under TeX and "off page" errors
under your driver.  

  You can dictate the edgescale for an individual figure by giving
the scale in brackets immediately after \SetLabels.  Thus, to
import into an article using say \Edge{100} a figure labelled using
another edgescale, say the original 1-by-1 default, you can use
\SetLabels[1]...\endSetLabels.

GETTING IT DOWN PAT

     Complicated labeling deserves the same respect as
complicated mathematics.  Do not expect it to come out perfect the
first time!  What is needed in either case is a mechanism to
repeatedly typeset troublesome pieces.

     One mechanism is always available.  One does complicated
labelling in a separate "test" file involving just the figure being
labelled;  a texpert will know how to \dump TeX's current state as
a temporary format that restarts rapidly at each retry.  Usually,
one then pastes the completed labelled figure back into the main
TeX file, but, of course, one can also \input it as an auxiliary
file.

     If you do not have a TeXpert at handy, here is a first
approximation to an efficient setup. By deletions reduce a copy
of your article to just a few lines before and after the figure.
Now label the figure, and finally, copy and paste the labelled
figure to the original article. Then copy the next figure to label
into this testbed and repeat. The TeXpert can improve the  speed
at which TeX starts up, by compiling a format specifically for
your article; just one caution: best NOT include in the format
ephemeral details of setup like \Set<mydriver>ArtSpecials (from
boxedeps.tex because this reads  figure dimensions which you may
change during your work session.

     An improved mechanism to repeatedly typeset troublesome
pieces is now available on the Macintosh; it is called LinoTeX;
see the same ftp sources.  It could be set up on many types
of computer.

     Before using labelfig.tex to attach labels to a graphics
object inserted using boxedeps.tex or BoxedArt.tex, make it a
firm rule to carefully adjust the bounding box using the trimming
commands of these packages, and also at least tentatively scale
and position the object. Beware of changing the grid inadvertently
after the labels have been positioned.  For example, correcting
the bounding box of a PostScript graphics object can foul up the
labels by changing the coordinate grid to which the labels are
attached. This is particularly true for the trimming  commands of
boxedeps.tex and BoxedArt.tex. However, as noted already, change
of scale is much less disruptive, and modest adjustments should be
well tolerated.

     Sometimes the labels protrude so far from the bounding box
of a figure that the figure has to be repositioned.  Best do this
by ad hoc spacing, say using \hglue and \vglue; altering the
bounding box would create a vicious circle.

     Remember that you are responsible for preventing labels
from overlapping. You are responsible for all label typography
including size and style. A label is really just about anything
that can be put in a TeX box. Note that spaces at the beginning
and end of labels will normally be suppressed; if you really want
them you must protect them with TeX braces.

     This package temporarily sets the \mathsurround parameter
of TeX to zero  while the labels are being affixed. This is done
because nonzero \mathsurround space would influence the position
of left and right aligned labels; then, when a texpert or printer
modifies mathsurround, diagram labeling might be disastrously
altered. There is a small price to pay involving labels that are
formatted as caption boxes including mathematics: you  may want or
need to specify an explicit mathsurround space within the caption
box; it will not influence anything outside.

     Those hostile to the use of * as separator between
the X and Y coordinates of label insertion points, are free to
impose another using \SetXYSeparator{<the new separator>}.  
Americans may prefer "," to "*" since they never use a 
comma as a decimal point; on the other hand, * may be more visible.

APPENDIX (I)  MERGING labelfig.tex LABELS INTO AN EPSF GRAPHICS OBJECT.

     As promised in the introduction, here is a recipe useful for
publishers. It works at least on Macintosh and at least for vectorized
graphics and Adobe type1 fonts.  (There is surely a similar recipe for
PCs under MSWindows.)

 (a)  Use boxedeps.tex utility to integrate the figure given by the eps
file, "x.eps" say, with a visible frame around it.  See
\ShowDisplacementBoxes command in boxedeps.tex.  To get precise results
automatically it is important to use the \Trim... commands of
boxedeps.tex making the "DisplacementBox" neatly fit the figure.

 (b)  Use the TeX printer driver and LaserWriter (versions >= 8.1.1) to
export to an EPSF the DVI page containing the integrated, labelled
figure. You now have an EPS file  "xx.eps"  that contains too much, and at
the wrong scale, and at wrong position.

 (c)  Convert the EPSF to an Adode Illustrator format EPSF using
the shareware utility called epsConvert by Sam Weiss
1993-- (currently $25).

 (d)  In Illustrator (or a compatible program), group the labels and the
"DisplacementBox"; copy them to the clipboard and paste them into "x.ps".
This step requires that all the label fonts be "visible to the Macintosh.

 (e)  Translate and scale the pasted group consisting of the labels plus
the "DisplacementBox" so as to make the "DisplacementBox" the bounding
box of (labelless) figure represented by "x.eps".  At this point the
labels will be correctly placed on the figure "x.eps".

 (f)  Ungroup and delete the "DisplacementBox".  The result is the
desired single EPS file, "x+.eps" say, It contains the original figure
plus its labels.  

     Using grouping and ungrouping appropriately in "x+.eps", a
publisher's staff can very efficiently improve label positions etc.

APPENDIX II)  SOME EXOTIC APPLICATIONS

     The grid of labelfig.tex is analogous to a light-table in
classical page makeup with wax or latex glue.  In principle, you
can use it to compose any page from its indivisible parts.  This
even has some of the artisanal charm of classical paste-up
provided you have a fast screen preview to make the process
"interactive".

     In practice labelfig.tex is a tool for nonstandard jobs.
Here are a few going beyond the labelling already discussed.

(I)  GRAPHICS INTEGRATION.

     This is accomplished by treating the imported graphics
objects as labels.  The underlying graphics object is then
typically an empty  \vbox to <dimension>{\vfill} in a TeX
\midinsert...\endinsert construction.  A label line
might be of the form

   (.1*.1) \\

The exact form of the special command varies from driver to
driver.  However, in the case of encapsulated PostScript graphics
(EPSF norm), by relying on boxedeps.tex, one can have the
following standard syntax (independant of driver  (see
boxedeps.doc for details.
  
  (.1*.1) \BoxedEPSF{MyFigure scaled <scale in mils>}\\

This may be slow since it requires TeX to read the PostScript
file to read bounding box using many complex macros.  So you
may want to try

  (.1*.1) \EPSFSpecial{MyFigure}{<scale in mils>}\\

which is fast and driver independant, but it squashes the
bounding box, normally to its lower left corner.

     Similarly for graphics of the Macintosh PICT norm ---
using BoxedArt.tex (same sources) in place of boxedeps.tex.

     This approach to integration is to be recommended when
one is assembling a composite graphics object.

 (II)  COMMUTATIVE DIAGRAM ENHANCEMENT

     Commutative diagrams or arrays of mathematical objects
connected by arrows of various sorts are common in mathematics.
The mathematical objects require the use of TeX.  Recently TeX
acquired a good collection of arrows of all slopes --- that of
LamSTeX --- plus pwerful macros to build the diagrams.

     However, even the LamSTeX collection is often
inadequate; it lacks for example double shafted arrows, dotted
arrows and curved arrows. Fortunately it is possible to produce
such arrows on an individual basis using sophisticated graphics
programs such as Illustrator and AldusFreehand (both serving
the EPSF norm) or using Metafont (with its public domain norm).
Since the creation of each new arrow is a work of love, you
probably want to limit the number of arrows by using LamSTeX
for most arrows. The 40K commutative diagram module of LamSTeX
has been adapted to work with AmSTeX and a copy may be posted
with LabelFig and related files. Unfortunately no one has yet
offered a version that works with Plain TeX or LaTeX.

       Suffice it here to say that when the exotic arrow has
been somehow imported into TeX, labelfig.tex treats it as a
label that one affixes to the commutative diagram.  Two other
steps will be treated in separate notes, namely the matter of
extracting the dimension specifications for the arrow and the
construction of the arrow --- for these steps are far from
unique and often depend intimately on your computer environment. 
Notes for the Macintosh-Textures-Illustrator combination are
found in the file ExoticArrows.doc.

 (III) NESTING 

Ingenuity pays off in exploiting labelfig.tex. One can
mix graphics and typography quite freely.  labelfig.tex is good
for freeform or overlapping arrangements, while boxedeps.tex (or
BoxedArt.tex) is best for regimented non-overlapping
arrangements --- and the two can be combined.

     The default behavior of labelfig.tex is not ideal 
for nesting objects, because to prevent trouble for beginners
the register for labels is globally cleared when \AffixLabels
concludes.  But there are switches available

      \LabelsGlobal      \LabelsLocal

which change this.  To understand this, extend the above test 
by something like:

 \LabelsLocal
 \SetLabels
    (.5*.5) AAA\\
 \endSetLabels

 {%%% Watch for influence of braces!!
 \SetLabels
    (.5*.5) ZZZ\\
 \endSetLabels
   \AffixLabels{\FirstQuadrant}
 }

   \AffixLabels{\FirstQuadrant}

     There are however potential pitfalls.  Neither
labelfig.tex nor boxedeps.tex has been tested under extreme
conditions. Problems may occur if their procedures are
indiscriminately nested. For boxedeps.tex (not labelfig.tex)
there is a precise cause for worry, namely many of its
variables are "global", which means that TeX braces will not
provide the protection one might expect.

COMMAND SUMMARY FOR labelfig.tex

  Here [...] means optional (one or zero)
       [...]* means any number of such constructs

  \SetLabels
    [[<P>](<X><Sep><Y>) <label> \\]*
  \endSetLabels
  \ShowGrid  % this makes the grid visible (once)
  \AffixLabels{<the figure>}

   --- <P> is tack position, one of eleven or empty
              order irrelevant

                   \L\T      \T      \R\T

                   \L\E      \E      \R\E

                     \L               \R

                   \L\B      \B      \R\B

   --- (<X><Sep><Y>) insertion point;
  <Sep> is separator, = * by default;
  \SetXYSeparator{<Sep>} changes it.
   <X> and <Y> are real numbers

  --- <label> a label to attach 

  --- <the figure> the figure to label 

  \GlobalLabels (default)     
  \LocalLabels  setting for nested constructs.

 \Grids makes ALL grids appear; \HideGrid then makes just next disappear.
 \noGrids returns to default.  The commands are always global.

 \GridLineWidth{<dimension>} adjusts width of grid lines. Default is very
small, to give "hairline" effect. If your grid lines are missing try
setting \GridLineWidth{1pt}.

 \Edges#1 globally changes the edge size of all grids to the numerical 
value #1, which must be 1, 10, 100, or 1000.  The default is 1.

VERSION HISTORY.
 --- Jan 1993: \Edges#1 and [??] option after \SetLabels
 --- July 1992: \Grids, \noGrids, \HideGrid;
       Gridlines become hairlines; \GridLineWidth{<dimension>}.
 --- Oct 1991, Jan 1992: \SetXYSeparator{<Sep>},  \LabelsGlobal,
       \LabelsLocal.
 --- July 1991: first release

Address for bugs and other feedback:

        Raymond S\'eroul
        IREM and Lab. de Typographie Informatise
        Univ. Rene Descartes
        Strasbourg

    Tel 33-88-41-63-45
    Email:  A18645@FRCCSC21.BITNET

        Laurent Siebenmann
        Mathematique, Bat. 425,
        Univ de Paris-Sud,
        91405-Orsay,
        France

    Tel 33-1-6941-7949; 
    Email: lcs@topo.math.u-psud.fr

\newtheorem{tetel}{Theorem}[section]
\newtheorem{defi}[tetel]{Definition}
\newtheorem{prop}[tetel]{Proposition}
\newtheorem{lemma}[tetel]{Lemma}
\newtheorem{rem}[tetel]{Remark}
\newtheorem{cor}[tetel]{Corollary}
\newtheorem{exa}[tetel]{Example}
\newtheorem{quest}[tetel]{Question}
\theoremstyle{definition}\newtheorem{exa1}[tetel]{Example}

\begin{document}

\newcommand{\diam}{\mathop{\textrm{diam}}}
\newcommand{\dist}{{\mathop{\textrm{dist}}}}
\newcommand{\lra}{\leftrightarrow}
\newcommand{\xlra}{\xleftrightarrow}
\newcommand{\xnlra}{\xnleftrightarrow}
\newcommand{\pc}{{p_c}}
\newcommand{\pt}{{p_T}}
\newcommand{\ptk}{{p_T^\mathsf{a}}}
\newcommand{\pl}{{\tilde{p}_c}}
\newcommand{\pe}{\tilde{p}_c^\mathsf{a}}
\newcommand{\pr}{\mathrm{\mathbb{P}}}
\newcommand{\pp}{\mu}
\newcommand{\ex}{\mathrm{\mathbb{E}}}
\newcommand{\ee}{\mathrm{\overline{\mathbb{E}}}}
\newcommand{\C}{\mathcal{C}}
\newcommand{\F}{\mathcal{F}}
\newcommand{\A}{\mathcal{A}}
\newcommand{\om}{{\omega}}
\newcommand{\ebd}{\partial_E}
\newcommand{\ivbd}{\partial_V^\mathrm{in}}
\newcommand{\ovbd}{\partial_V^\mathrm{out}}
\newcommand{\q}{q}
\newcommand{\Z}{\mathbb{Z}}
\newcommand{\T}{\mathbb{T}}
\newcommand{\TT}{\mathcal{T}}
\newcommand{\GG}{\mathbb{G}}
\newcommand{\eps}{\varepsilon}
\newcommand{\fss}{\mathcal{S}}

\title{On percolation critical probabilities and unimodular random graphs}

\author{Dorottya Beringer
\and 
G\'abor Pete
\and
\'Ad\'am Tim\'ar}
\date{}
\maketitle

\begin{abstract}
We investigate generalisations of the classical percolation critical  probabilities $\pc$, $\pt$ 
and the critical probability $\pl$ defined by Duminil-Copin and Tassion \cite{DCT} 
to bounded degree unimodular random graphs. 
We further examine Schramm's conjecture in the case of unimodular random graphs: 
does $\pc(G_n)$ converge to $\pc(G)$ if $G_n\to G$ in the local weak sense? 
Among our results are the following: 
\begin{itemize}
\item $\pc=\pl$ holds for bounded degree unimodular graphs. However, there are unimodular graphs with sub-exponential volume growth and $\pt<\pc$; i.e., the classical sharpness of phase transition does not hold.
\item We give conditions which imply $\lim\pc(G_n)= \pc(\lim G_n)$. 
\item There are sequences of unimodular graphs such that $G_n\to G$ but $\pc(G)>\lim \pc(G_n)$ or $\pc(G)<\lim \pc(G_n)<1$.
\end{itemize}
As a corollary to our positive results, we show that for any transitive graph with sub-exponential volume growth there is a sequence $\TT_n$ of large girth bi-Lipschitz invariant subgraphs such that $\pc(\TT_n)\to 1$. It remains open whether this holds whenever the transitive graph has cost 1.
\end{abstract}

\section{Introduction}\label{s.intro}

\subsection{Motivation and results}\label{ss.results}

There are several definitions of the critical probability for percolation on the lattices $\Z^d$, which have turned out to be equivalent not only on $\Z^d$, but also in the more general context of arbitrary transitive graphs \cite{M, AB, Gr, AV, DCT, DCTZd}. One of our goals is to investigate the relationship between these different definitions when the graph $G$ is an ergodic unimodular random graph \cite{BS, AL07}, which is the natural extension of transitivity to the disordered setting.
We examine the generalisations of $\pc=\sup\{p: \pr_p(\textrm{there is an infinite cluster})=0\}$, 
$\pt=\sup\left\{p: \ex_p(|\C_o|)<\infty\right\}$ and 
$\pl$, defined by Duminil-Copin and Tassion \cite{DCT}. The last quantity was in fact designed to give a simple new proof of $p_c=p_T$ for transitive graphs, and to address the question of locality of critical percolation: whether the value of $p_c$ depends only on the local structure of the graph.

More precisely, Schramm's ``locality conjecture'', stated first explicitly in \cite{BNP}, says that $\pc(G_n)\to \pc(G)$ holds
whenever $G_n$ is a sequence of vertex-transitive infinite graphs such that $G_n$ converges locally to $G$ (i.e., for every radius $r$, the $r$-ball in $G_n$, for $n$ large enough, is isomorphic to the $r$-ball in $G$) and $\sup_n\pc(G_n)<1$. Typically, however, the natural setting for such locality statements is not the class of transitive graphs, but the class of unimodular random graphs. Indeed, there are several interesting probabilistic quantities, most often related in some way to random walks, which have turned out to possess locality, mostly in the generality of unimodular random graphs: see \cite{BS, Lyo05, LN, CsF, BSzV, GL} for specific examples, and \cite[Chapter 14]{PGG} for a partial overview.  Therefore, it is natural to investigate Schramm's conjecture in the setup of unimodular random graphs and see what the proper notion of critical probability may be from the point of view of locality.

The conjecture has been proved for some special transitive graphs.
Grimmett and Marstrand \cite{GM} proved that
$\pc\left( \Z^2\times \{-n,\dots ,n\}^{d-2} \right)\xrightarrow{n\to\infty} \pc(\Z^d)$.
Benjamini, Nachmias and Peres \cite{BNP} 
verified that the convergence holds if $(G_n)$ is a sequence of $d$-regular graphs with large girth and
Cheeger constants uniformly bounded away from 0.
Martineau and Tassion \cite{MT} proved that the convergence holds if $(G_n)$ is a sequence of Cayley graphs of Abelian groups converging to a Cayley graph $G$ of an Abelian group, and $\pc(G_n)<1$ for all $n$. The inequality
\[ \liminf_{n\to\infty}\pc(G_n)\geq \pc(G) \]
is known for any convergent sequence of transitive graphs; see \cite[Section 14.2]{PGG}, and \cite{DCT}. 
Given the scarcity of transitive examples, it is a natural wish to try and find classes of unimodular graphs that satisfy the locality or at least the lower semicontinuity of the critical probability. 

In Subsection~\ref{ss.pcdef}, we define the generalized critical probabilities
$\pc$, $\pt$, $\pl$, $\ptk$, and $\pe$ for unimodular random graphs; somewhat
simplistically saying, the first three will be quenched versions of the
quantities mentioned above, while the last two will be annealed versions.  

In Section~\ref{s.pcrel}, we examine the relationship between these different generalizations. 
The main positive result of this section, used many times in the rest of the paper, is the following: 

\begin{tetel}\label{pl=pc}
If $(G,o)$ is a bounded degree unimodular random rooted graph, then $\pc(G)=\pl(G)$ holds a.s., where $\pl$ is the quantity introduced by Duminil-Copin and Tassion \cite{DCT}; see~(\ref{pltrans}) below.
\end{tetel}

%We discuss the motivation and the details of the definition of $\pl$ in Section~\ref{ss.pcdef}. 
Our further results on the relationship of the different definitions of critical probabilities are summarized in Table~\ref{tab}. The one sentence summary is that although $\pc=\pl$ always holds, otherwise almost anything can happen, unless the random graph satisfies some very strong uniformity conditions; one that we call ``uniformly good'' suffices for most purposes. 
The notion of uniformly good unimodular graphs (see Definition \ref{unifgood}) captures the property of the original definition of $\pl$ that there is a bounded size witness for $p$ being less than $\pl$. 
This class of graphs includes all quasi-transitive unimodular graphs and unimodular trees of sub-exponential growth. %\gabor{only trees, not subexp graphs, right?}

In Section~\ref{s.loc} we investigate the extension of Schramm's conjecture to unimodular random graphs: 
does $\pc(G_n)$ converge to $\pc(G)$ if $G_n\to G$ in the local weak sense
(i.e., the laws of the $r$-balls in $G_n$ converge weakly, for every $r$) and
$\sup \pc(G_n)<1$? First we note that locality holds for
unimodular Galton--Watson trees with bounded degrees, but not in general. 

\begin{exa}
Let $X_n$ and $X$ be uniformly bounded non-negative integer valued random variables with mean larger than 1. 
Denote by $UGW_\infty(X_n)$ and $UGW_\infty(X)$ the unimodular Galton--Watson trees with these offspring distributions, conditioned to be infinite. If $UGW_\infty(X_n) \to UGW_\infty(X)$ in the local weak sense, then $\pc(UGW_\infty(X_n))\to \pc(UGW_\infty(X))$. 
\end{exa}

We discuss this family of graphs in more detail in Example~\ref{ugwt}. This example motivates our investigations on the locality of the critical probability in the class of unimodular random graphs and it shows that it is natural to restrict one's attention to bounded degree graphs. 

In Subsection~\ref{ss.semi}, we prove some general positive results: lower semicontinuity of $\pe$ and the following two propositions, giving some particular settings where locality of $\pc$ holds:

\begin{prop}\label{liminf}
Let $G$ be a uniformly good unimodular random graph. Furthermore, let $G_n$ be uniformly bounded degree unimodular random graphs converging to $G$ in the local weak sense, in a \emph{uniformly sparse} way: there is a positive integer $k$ such that for each $n$ there is a coupling $\nu_n$ of $\mu_G$ and $\mu_{G_n}$ such that $G\subseteq G_n$ and 
there is a sequence of positive integers $r_n\to\infty$ that satisfies 
$|(E(G_n)\setminus E(G))\cap B_{G_n}(o,r_n)|\leq k$ $\nu_n$-almost surely. Then
$\lim_{n\to\infty}\pc(G_n)= \pc(G).$
\end{prop}

Although uniformly good unimodular graphs are not much more general than quasi-transitive graphs, the main point of this proposition is that it gives examples satisfying locality, beyond transitive graphs and unimodular Galton--Watson trees; see, e.g., Example~\ref{pmexa}. Also, it draws attention to how fragile locality is in the realm of unimodular random graphs: in Subsection~\ref{ss.loc.countex}, we show by examples that neither uniformly sparse convergence, nor a uniformly good limit suffices alone for locality: there are such sequences of unimodular random graphs with $G_n\to G$ but $\pc(G)>\lim \pc(G_n)$ or $\pc(G)<\lim \pc(G_n)<1$. 

In the quite special setting of unimodular trees of uniform subexponential growth (see Definition~\ref{subexp.def}), the assumption of uniformly sparse convergence from Proposition~\ref{liminf} can be relaxed:

\begin{prop}\label{subexp}
If $G$ is a bounded degree unimodular random tree with uniformly subexponential volume growth, then all five critical percolation densities equal 1, and $G$ is uniformly good. If $G_n$ is a sequence of bounded degree unimodular random graphs with uniformly subexponential volume growth and girth tending to infinity, then $\pc(G_n)$, $\pl(G_n)$, $\pe(G_n)$ all tend to $1$.
\end{prop}

A corollary to this result is that if $G$ is a transitive graph of subexponential volume growth, then there exists a sequence of invariant bi-Lipschitz spanning subgraphs $G_n$ such that $\pc(G_n) \to 1$. As we will explain in Section~\ref{s.cost}, this is a strengthening of the simple fact that groups of subexponential growth have cost 1, as defined in \cite{Lcost}, studied further in \cite{Gcost,GICM}. We do not know if this strengthening holds for all groups of cost 1, which class includes, besides all amenable groups, direct products $\GG\times\Z$ for any group $\GG$, and $SL(d,\Z)$ with $d\ge 3$. A related question is whether every amenable transitive graph has an invariant random Hamiltonian path. This is the invariant infinite version of what is known as Lov\'asz' conjecture, namely, that every finite transitive graph has a Hamiltonian path, even though he has not conjectured a positive answer. The best general results seem to be \cite{B} and \cite{PR}.
\medskip

Our positive results notwithstanding, a key conclusion of our work seems to lie in the counterexamples: there appears to be no perfect definition of a ``critical density'' that would make locality a robust phenomenon, true for a large class of unimodular random graphs and thus possibly more accessible for a proof in the transitive case. %\gabor{Ez jo otlet? Lasd meg ref1 reply.}

\subsection{Notation}\label{ss.not}

\noindent{\bf Graphs.} We always consider locally finite and rooted graphs. The root is denoted by $o$.
We denote by $e^-$ and $e^+$ the endpoints of the (directed) edge $e$. 
When a subgraph $S$ is given (maybe implicitly) and it contains exactly one endpoint of $e$, 
then we denote that endpoint by $e^-$. 
We write $x\sim y$ if $x$ and $y$ are adjacent vertices in $G$.
We will use $\dist_G(x,y)$ for the graph distance between the vertices $x$ and $y$ in the graph $G$.
We denote by $B_G(o,r)$ the ball around $o$ of radius $r$ in $G$, i.e.,
the subgraph induced by the vertex set $\{x\in V(G): \dist_G(o,x)\leq r\}$. 
For any subset $S$ of the vertices, let 
$\ebd S:=\left\{ e\in E(G): e^-\in S, e^+\notin S \right\}$ be the \emph{edge boundary} of $S$, let
$\ivbd S:=\left\{ x\in S: \exists y\sim x, y\notin S \right\}$  be the \emph{internal vertex boundary} of $S$, and let
$\ovbd S:=\left\{ x\notin S: \exists y\sim x, y\in S \right\}$ be the
\emph{outer vertex boundary} of $S$.
For a rooted graph $(G,o)$ we denote by $\fss(G)$ the set of finite subsets of $G$ which
contain the root $o$. 

Several of our examples will use percolation on $\Z^2$. The subgraph spanned by the box $[-n,n]^2$ will be denoted by $Q_n$. We will also use the standard dual percolation on the dual lattice $(\Z+\frac12)^2$.

When we talk about invariant random subgraphs of a Cayley graph $\Gamma$ of a group $\GG$, we will always mean that the measure on subgraphs is invariant under the natural action of $\GG$. When we talk about invariant random subgraphs of a transitive graph $\Gamma$, with no group action specified, then we mean invariance under the automorphism group $Aut(\Gamma)$.

\medskip
\noindent{\bf Unimodular random graphs.} Let $\mathcal{G}_{\star}$ be the space of isomorphism classes of locally finite labeled rooted graphs, and let $\mathcal{G}_{\star\star}$ be the space of isomorphism classes of locally finite labeled graphs with an ordered pair of distinguished vertices, each equipped with the natural local topology: two (doubly) rooted graphs are ``close'' if they agree in ``large'' neighborhoods of the root(s). 
If $(G,o)$ is a random rooted graph, then denote by $\pp_G$ the distribution of it on $\mathcal{G}_{\star}$, and 
let $\ee_G$ be the expectation with respect to $\mu_G$. 
We omit the index $G$ from this notation if it is clear what the measure is. 

\begin{defi}[\cite{AL07}, Definition 2.1]
We say that a random rooted graph $(G,o)$ is \emph{unimodular} if it obeys the Mass Transport Principle:
\[ \ee_G\left(\sum_{x\in V(\om)}f(\om,o,x)\right) = \ee_G\left(\sum_{x\in V(\om)}f(\om,x,o)\right) \]
for each Borel function 
$f: \mathcal{G}_{\star\star}\to [0,\infty]$.
\end{defi}

There are several other equivalent definitions; see \cite[Definition 14.1]{PGG}. Also, it is an open question if this class is strictly larger than the class of sofic measures: the closure of the set of finite graphs under local weak convergence.

An important class of unimodular graphs consists of Cayley graphs of finitely generated groups 
and of invariant random subgraphs of a Cayley graph:

\begin{prop}[\cite{AL07}, Remark 3.3]\label{caysub}
Let $\Gamma$ be a Cayley graph of a finitely generated group and let $o$ be a vertex of $\Gamma$.
If $G$ is a random subgraph of $\Gamma$ that is invariant under the action of the group, 
then $(G,o)$ is unimodular.
\end{prop}

The class of unimodular probability measures is convex. 
A unimodular probability measure is called \emph{extremal} if 
it cannot be written as a convex combination of other unimodular probability measures. 

\medskip
\noindent{\bf Percolation.} For simplicity, we will consider only {\it bond} percolation processes on unimodular random graphs. For a fixed instance $\om$ of the random graph $G$ let $\pr_p^\om$ be the probability measure 
obtained by the Bernoulli($p$) bond percolation on $\om$ and  let $\ex_p^\om$ be the expectation with respect to $\pr_p^\om$. The percolation \emph{cluster} (i.e., the connected component) of the root $o$ will be $\C_o$.

\subsection{Critical probabilities}\label{ss.pcdef}

The long studied critical probabilities $\pc=\sup\left\{p: \pr_p(|\C_o|=\infty)=0\right\}$ 
first defined by Hammersley 
and $\pt=\sup\left\{p: \ex_p(|\C_o|)<\infty\right\}$ introduced by Temperley 
have natural generalizations to extremal unimodular random graphs. 
Let $(G,o)$ be an extremal unimodular random graph. 
In this case the critical probability $\pc(\om)$ of an instance of $(G,o)$ is almost surely a constant 
and the same holds for $\pt$ (see \cite{AL07}, Section 6.). 
Hence one can define 
\begin{align*}
\pc&=\inf\left\{p: \pp\left(\pr_p^\om\left(|\C_o|=\infty\right)>0\right)=1\right\} \\
 &=\sup\left\{p: \pp\left(\pr_p^\om\left(|\C_o|=\infty\right)=0\right)=1\right\}
\end{align*}
and
\begin{equation*}
\begin{aligned}
\pt&=\sup\left\{p: \pp\left(\ex_p^\om\left(|\C_o|\right)<\infty\right)=1\right\}\\
 &=\inf\left\{p: \pp\left(\ex_p^\om\left(|\C_o|\right)=\infty\right)=1\right\}.
\end{aligned}
\end{equation*}
It may happen that although $\ex_p^\om\left(|\C_o|\right)<\infty$ for $\mu$-almost every $\om$,  
the expectation of these quantities with respect to $\mu$ is infinite. 
This provides a second natural extension of $\pt$ to unimodular random graphs defined using the average size of $\C_o$: 
\begin{align*}
\ptk&=\sup\left\{p: \ee\left(\ex_p^\om\left(|\C_o|\right)\right)<\infty \right\} \\
 &= \inf\left\{p: \ee\left(\ex_p^\om\left(|\C_o|\right)\right)=\infty\right\}.
\end{align*}
It follows from the definitions that $\pc\geq \pt\geq \ptk$. 
It is known that $\pc=\pt$ in the case of transitive graphs; see \cite{M, AB, AV, DCT}. 
For unimodular random graphs (even with sub-exponential volume growth), the three critical probabilities can differ; 
we will present such graphs in Examples \ref{subexpcanopy} and \ref{ptk}. 

Duminil-Copin and Tassion \cite{DCT} introduced the following local quantity for transitive graphs:
let $(G,o)$ be a rooted graph, $S\in\fss(G)$ be a finite subgraph containing the root, and define
\[\phi_p(S):=\sum_{e\in\ebd S} p \, \pr_p(o\xlra{S}e^-)\,,\]
the expected number of open edges on the boundary of $S$
such that there is an open path from $o$ to $e^-$ in $S$.
Then, they defined the critical probability
\begin{equation}\label{pltrans}
\begin{aligned}
\pl&:=\sup\{p: \textrm{there is an } S\in\fss(G) \textrm{ s.t. } \phi_p(S)<1\}\\
 &=\inf \{p:  \phi_p(S)\geq 1 \textrm{ for all } S\in\fss(G)\}\,. 
\end{aligned}
\end{equation}
They proved that transitive graphs satisfy $\pc=\pl$. 

How to generalize this definition to unimodular random graphs is not a priori clear. The simplest 
way to define a similar critical probability seems to be a quenched version: find a suitable $S_\om\in\fss(\om)$ for almost every configuration $\omega$. For a subgraph $S\in\fss(\om)$ denote by  
\begin{equation}
\label{phi} \phi_p^\om(S):=\sum_{e\in\ebd S} p\,\pr_p^\om\left(o\xlra[S]{\om,p}e^-\right)  
\end{equation}
the expected number of open edges on the boundary of $S$ in $\om$ 
such that there is an open path from $o$ to $e^-$ in the percolation on $\om$ with parameter $p$.
Then let
\begin{equation}\label{plque}
\pl:=\sup\left\{p:
\pp\left(\left\{\omega : \exists S_\om\in\fss(\om) \textrm{ s.t. } \phi_p^\om(S_\om)<1 \right\}\right)=1\right\}. 
\end{equation}

\begin{rem}\label{Somx}
Suppose $p$ satisfies the following: for almost every $\om$ there is an $S_\om\in\fss(\om)$ with $\phi_p^\om(S_\om)<c$. 
Then unimodularity implies \cite[Lemma 2.3.]{AL07} that
for almost every $\om$ and every vertex $x$ there is some finite connected set $S_{\om,x} \ni x$
such that
\[\phi_p^{\om,x}(S_{\om,x}):=p\sum_{e\in\ebd S_{\om,x}}\pr\left(x\xlra[S_{\om,x}]{\om,p} e^-\right)<c.\]
\end{rem}

In the original definition~\eqref{pltrans} of $\pl$, there is no control on what the set $S$ could be, 
which makes the definition rather ineffective. This becomes particularly problematic in the random graph case~\eqref{plque}, where a bad neighborhood of $o$ may force $S_\omega$ to be huge and hard to find. However, it will follow from our Lemma~\ref{unifgoodexp} that, for transitive graphs, the existence of an $S$ with $\phi_p(S)<1$ is
equivalent to the existence of a positive integer $r$ with $\phi_p(B(o,r))<1$. 
This provides a second natural extension of the definition of $\pl$ to the random case: 
we consider the ball of radius $r$ in the random graph $\om$ and we take the expectation of $\phi_p^\om(B_\om(o,r))$ with respect to $\mu$. 
Then the following critical probability is another extension of the definition of $\pl$: 
\begin{equation}\label{plann}
\pe:=\sup\{p: \exists r \textrm{ such that }\ee\left(\phi_p^\om(B_\om(o,r))\right)<1 \}.
\end{equation}

\subsection{Operations preserving unimodularity}\label{ss.opera}

We give now the general description of some operations
on the space of unimodular graphs that we will use in our counterexamples in Subsections~\ref{ss.pc.countex} and~\ref{ss.loc.countex}. This subsection is not necessary to understand the positive results of the paper. 

Some of our examples arise from Cayley graphs using operations from $\mathcal{G}_\star$ to $\mathcal{G}_\star$. 
One of these operations is the \emph{edge replacement} defined in \cite{AL07}, Example 9.8: 
we replace each edge of a unimodular graph $G$ by a finite graph with two distinguished vertices corresponding to the endpoints of the edge, then we find the correct new distribution for the root that makes the measure unimodular. If the finite graphs are random, each must have finite expected vertex size. 
In this section, we define further operations, called \emph{vertex replacement} and \emph{contraction}, and 
we prove that if the initial graph is a unimodular labeled graph with appropriate labels, 
then the resulting graph by such an operation is also unimodular.
\medskip

\noindent\textbf{Vertex replacement.}
Let $(\Gamma,o)$ be a unimodular random labeled graph with distribution $\mu$, where the labels are in the form $(G_x, \varphi_x)$, where 
$G_x$ is a finite graph and $\varphi_x$ is a map from $\{(x,y)\in E(\om): y\sim x\}$ to $V(G_x)$. 
If the labeling satisfies 
$\ee_{\mu} |V(G_o)|<\infty$, 
then we can define the following rooted random graph $H(\Gamma)$: 
we choose $(\Gamma, o, \{(G_x,\varphi_x):x\in V(\Gamma)\})$ with respect to the probability measure $\mu$ biased by $|V(G_o)|$, and replace each vertex $x$ of $\Gamma$ by the graph $G_x$ and each edge $e$ of $\Gamma$ 
by the edge $\{\varphi_{e^-}(e),\varphi_{e^+}(e)\}$. 
Let the root $o'$ of $H(\Gamma)$ be a uniform random vertex of $V(G_o)$. 
Denote the law of $(H(\Gamma), o')$ by $\mu'$.  

We claim that if $\mu$ is unimodular with $\ee_{\mu} |V(G_o)|<\infty$, then $\mu'$ is also unimodular. 
Let $f(\om, u, v)$ be a Borel function from $\mathcal{G}_{\star\star}$ to $[0,\infty]$ 
and let 
\[\bar{f}(\bar{\om}, x, y):=\frac{1}{\ee_{\mu}|V(G_o)|} \sum_{u\in V(G_x), v\in V(G_y)} f(H(\bar{\om}), u, v) \]
which is an isomorphism-invariant Borel function on the subspace of $\mathcal{G}_{\star\star}$ that consists of graphs with labels of the above form. 
We show that $\mu'$ obeys the Mass Transport Principle: 
\begin{align*}
\int \sum_{v\in V(\om)} f(\om, o', v) d\mu'(\om,o') 
= &\int \sum_{o'\in V(G_o), v\in V(H(\bar{\om}))} 
\frac{1}{|V(G_o)|}f(H(\bar{\om}), o', v) \frac{|V(G_o)|}{\ee_{\mu}|V(G_o)|} d\mu(\bar{\om},o) \\
= &\int \sum_{x\in V(\bar{\om})}\sum_{o'\in V(G_o), v\in V(G_x)} 
\frac{1}{\ee_{\mu}|V(G_o)|}f(H(\bar{\om}), o', v) d\mu(\bar{\om},o) \\
= &\int \sum_{x\in V(\bar{\om})} \bar{f}(\bar{\om},o,x) d\mu(\bar{\om},o)  \\
= &\int \sum_{x\in V(\bar{\om})} \bar{f}(\bar{\om},x,o) d\mu(\bar{\om},o)  \\
= &\int \sum_{v\in V(\om)} f(\om, v, o') d\mu'(\om,o'). 
\end{align*}
\medskip

\noindent\textbf{Contraction.} 
Let $(\Gamma,o)$ be a unimodular random edge-labeled graph with distribution $\mu$, where the labels of the edges are 0 or 1. 
We denote by $G$ the random subgraph of $\Gamma$ spanned by all the vertices and the edges with label 1. 
For a vertex $x$ of $\Gamma$ let $\C_x$ be the connected component of $x$ in $G$. 
We define the contracted graph $H(\Gamma)$: 
in practice, this is what we get by identifying every vertices in the same component of $G$. 
More formally, first we choose $(\Gamma, o, G)$ with respect to the distribution $\mu$ biased by $\frac{1}{|\C_o|}$. 
The vertices of $H(\Gamma)$ are the connected components of $G$ and 
we join two vertices by an edge iff there is an edge in $\Gamma$ which connects the two components. 
Let the root $o'$ of $H(\Gamma)$ be the connected component $\C_{o}$.  
Denote the law of $(H(\Gamma), o')$ by $\mu'$. 

We claim that if $\mu$ is unimodular then $\mu'$ is also unimodular. 
Let $f(\om, u, v)$ be a Borel function from $\mathcal{G}_{\star\star}$ to $[0,\infty]$ 
and let 
\[\bar{f}(\bar{\om}, x, y):=\frac{1}{|\C_x||\C_y|} f(H(\bar{\om}), \C_x, \C_y) \]
which is an isomorphism-invariant Borel function on the subspace of $\mathcal{G}_{\star\star}$ 
that consists of graphs with edges labeled by 0 or 1, 
such that the subgraph defined by the edges with label 1 consists of finite components. 
We show that $\mu'$ obeys the Mass Transport Principle: 
\begin{align*}
\int \sum_{v\in V(\om)} f(\om, o', v) d\mu'(\om,o')
=& \int \sum_{x\in V(\bar{\om})} 
\frac{1}{|\C_x|}f(H(\bar{\om}), \C_o, \C_x) \frac{1}{|\C_o|\ee_{\mu}\left(\frac{1}{|\C_o|}\right)} d\mu(\bar{\om},o)  \\
=& \frac{1}{\ee_{\mu}\left(\frac{1}{|\C_o|}\right)} \int \sum_{x\in V(\bar{\om})} 
 \bar{f}(\bar{\om}, o, x) d\mu(\bar{\om},o)  \\
=& \frac{1}{\ee_{\mu}\left(\frac{1}{|\C_o|}\right)} \int \sum_{x\in V(\bar{\om})} 
 \bar{f}(\bar{\om},x, o) d\mu(\bar{\om},o)  \\
=& \int \sum_{v\in V(\om)} f(\om,v, o') d\mu'(\om,o'). 
\end{align*}

\section{Relationship of the critical probabilities}\label{s.pcrel}

We will start by proving Theorem \ref{pl=pc} that states that all bounded degree unimodular graphs satisfy $\pc=\pl$. This theorem will be used in many of our further results.  

In the transitive case, the quantity $\phi_p(S)$ in the definition of $\pl$ can be used to give a short proof (see \cite{DCT}) 
of Menshikov's theorem \cite{M}: 
if $\Gamma$ is a transitive graph and $p<\pc(\Gamma)$, then there exist a $\varphi(p)$ such that 
\begin{equation}\label{RadExpDec}
\pr_p\left(o\lra B(o,r)^c\right)\leq e^{-\varphi(p)r}.
\end{equation} 
If a graph satisfies this exponential decay for each $p<\pc$ and has sub-exponential volume growth, then it is easy to see that $\pt=\pc$. 
In Lemma \ref{unifgoodexp}, we give a condition for unimodular random graphs that implies (\ref{RadExpDec}), and we prove in Corollary~\ref{pc=pt} that this condition implies $\pc=\pt=\ptk$ if the graph has uniform sub-exponential volume growth. However, in Examples~\ref{subexpcanopy} and~\ref{ptk} we present unimodular random graphs with uniform polynomial volume growth and $\pt< \pc$ and $\ptk<\pt$, respectively. This shows that Menshikov's theorem is not true in the generality of unimodular graphs. 

The results of this section are summarized in Table~\ref{tab}.

\begin{table}[h!]
\centering
\begin{tabular}{|c| l|}
\hline
$\pl=\pc$ & bounded degree \\
\hline 
$\pc\geq \pt \geq \ptk$ & always \\
$\pc=\ptk$ & bounded degree uniformly good with sub-exp.~growth \\
$\pc> \pt$ & Example \ref{subexpcanopy}, with polynomial growth \\
$\pt>\ptk$ & Example \ref{ptk}, with polynomial growth\\
\hline 
$\pc\leq \pe$ & bounded degree uniformly good \\
$\pc>\pe$ & Example \ref{pc>pe}, not uniformly good\\
$\pc <\pe$ & Example \ref{pc<pe}, uniformly good\\
\hline
\end{tabular}
\caption{Relationship of the critical probabilities}
\label{tab}
\end{table}

\subsection{Positive results}\label{ss.pc.pos}

Our first result is indispensable to the rest of the paper. 
The second part of the proof is a slight modification of the proof in \cite{DCT} for our setting, while the first part depends on new ideas. 
The main difficulty is that we cannot find isomorphic sets $S_{\om,x}$ for different vertices $x$, and hence 
we cannot bound $\pr_p\left(o\lra B(o,r)^c\right)$ in terms of $r$. 
We build instead a tree $T^\om$ using the sets $S_{\om,x}$, and bound the probability that the subtree given by the percolation survives. The survival of that subtree is equivalent to the infinite size of the cluster of the root in the percolation on $G$.  

\begin{pfo} {\it Theorem \ref{pl=pc}.}
We prove first that $\pl\leq \pc$. Fixing $p<\pl$, we will show that $p\leq \pc$. 
We claim that there exists a constant $c=c(p)<1$ such that we can find for almost every $\omega$ 
a set $S_\om\in\fss(\om)$ that satisfies $\phi_p^\om(S_\om)\leq c$.
Let $p':=\frac{p+\pl}{2}<\pl$.
Let $S_\om\in\fss(\om)$ be such that $\phi_{p'}^\om(S_\om)<1$.
The sets $S_\om$ satisfy
\begin{align*} 
\phi_{p}^\om(S_\om)=\sum_{e\in\ebd S_\om}p\pr_{p}^\om(o\lra e^-) &\leq
\frac{p}{p'}\sum_{e\in\ebd S_\om}p'\pr_{p'}^\om(o\lra e^-) \\
&= \frac{p}{p'}\, \phi_{p'}^\om(S_\om) \leq \frac{p}{p'}=:c\,.
\end{align*}
Recall the definition of $\phi_p^{\om,x}(S_{\om,x})$ from Remark \ref{Somx}.
Unimodularity implies that almost every $\om$ satisfies the following: 
for each $x\in\om$ there is a set $S_{\om,x}$ containing $x$ such that $\phi_p^{\om,x}(S_{\om,x})\leq c$. 
Fix such an $S_{\om,x}$ in an arbitrary measurable way.

Fix $\om$ and denote by $T^\om$ the following recursively defined tree:
the vertices of the tree are finite sequences of vertices of $\om$.
The root of the tree is $(o)$.
If $(x_0, x_1, \dots , x_k)$ is a vertex of $T^\om$,
its children are the sequences $(x_0, x_1, \dots , x_k, x_{k+1})$ 
such that for all $j=1, \dots,{k+1}$, we have $x_j\in\ovbd S_{\om,x_{j-1}}$, and 
there exist vertices $x'_j\in\ivbd S_{\om,x_{j-1}}$ such that $x'_j\sim x_j$, with paths from $x_{j-1}$ to $x'_j$ in $S_{\om, x_{j-1}}$
that are disjoint from each other and from the edges $\{x'_j, x_j\}$, as $j=1,\dots ,k+1$.
We say that the union of the above paths and edges is a \emph{good path} through $x_0, x_1, \dots , x_k, x_{k+1}$. See Figure~\ref{goodpath}.
Denote by $L_n:=\{ (x_0, x_1, \dots , x_n)\in T^\om \}$ the vertex set of $T^\om$ on the $n$th level. 

\begin{figure}[htbp]
\SetLabels
(0.3*0.05) $\omega$\\
(0.12*0.6) $x_0$\\
(0.05*0.95) \color{Maroon}{$S_{\omega,x_0}$}\\
(0.35*0.63) $x_1'$\\
(0.43*0.49) $x_1$\\
(0.45*1) \color{Magenta}{$S_{\omega,x_1}$}\\
(0.53*0.33) $x_2'$\\
(0.75*0.42) $x_2$\\
(0.95*0.8) \color{Aquamarine}{$S_{\omega,x_2}$}\\
\endSetLabels
\centerline{
\AffixLabels{
   \includegraphics[width=0.5\linewidth]{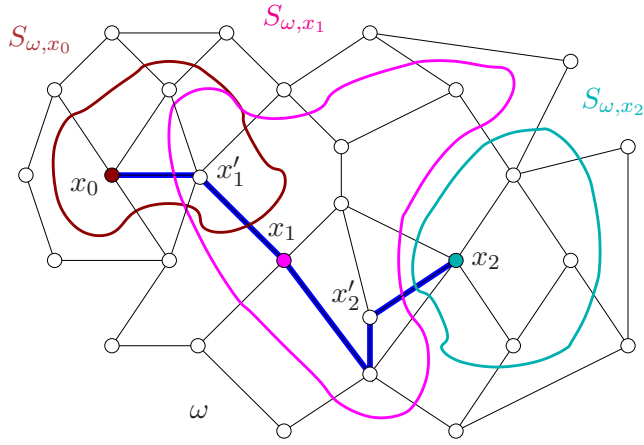} % requires the graphicx package
   }}
   \caption{A ``good path'' that gives the vertex $(x_0,x_1,x_2)$ of $T^\om$.}
   \label{goodpath}
\end{figure}

Let $T^\om(p)$ be the random subtree of $T^\om$ defined in a similar way
but allowing only good paths that are open in Bernoulli($p$) percolation on $\om$. 
%\gabor{gondolom, $S_{\om,x}$ helyett nem veszunk uj $S_{\om,p,x}$ halmazokat, szoval ez a megfogalmazas nem jo.}
It is easy to check that in fact $T^\om(p)\subseteq T^\om$. 
Denote by $L_n(p)$ the set of vertices of $T^\om(p)$ in the $n$th level. 
A self-avoiding infinite ray inside the $p$-percolation configuration gives rise to a growing sequence of good paths in the percolated $\om$, therefore if the cluster of the origin in the $p$-percolation on $\om$ is infinite, then there is an infinite path in $T^\om(p)$. Conversely, an infinite path in $T^\om(p)$ corresponds to an infinite growing sequence of open good paths in the $p$-percolated $\om$, which are necessarily parts of an infinite component containing the origin.
%Since any self-avoiding infinite ray inside the $p$-percolation configuration gives rise to an infinite good path that survives percolation \gabor{a megfogalmazas osszecsiszolando a $T^\om(p)$ definiciojaval}, and vice versa, the event that the cluster of the origin in the $p$-percolation on $\om$ is infinite coincides with the event that there exists an infinite path in the tree $T^\om(p)$ (that is, $T^\om(p)$ survives). 

We claim that for almost every $\om$ the expected number of vertices in $L_n(p)$ converges to 0 as $n\to \infty$.
More precisely, the expectation of the number of vertices in $L_n(p)$ decreases exponentially in $n$.
In the first two inequalities we use the notation $\square$ for the occurence of events on disjoint edge sets
and we apply the BK inequality (\cite{Gr}, Theorem 2.12). 
We denote the event $\Big\{x_0\xlra[B]{\om, p}x_k\textrm{ by a good path through }x_0, x_1, \dots ,x_k\Big\}$ 
by $\Big\{x_0\xlra[B,(x_0, x_1, \dots ,x_k)]{\om, p}x_k\Big\}$.
\begin{align*}
\ex^\om\big(|L_n(p)|\big)&
=\sum_{(x_0,\dots, x_n)\in L_{n}} \pr^\om\left( x_0\xlra[S,(x_0,\dots ,x_n)]{\om, p}x_n\right)\\
 &\leq \sum_{\substack{(x_0,\dots, x_{n-1})\in L_{n-1} \\ e\in\ebd S_{\om,x_{n-1}}}}
\pr^\om\left( \Big\{x_0\xlra[S,(x_0,\dots, x_{n-1})]{\om, p}x_{n-1}\Big\} \square \big\{e \textrm{ is open}\big\} \square 
\Big\{x_{n-1}\xlra[S_{\om,x_{n-1}}]{\om, p} e^-\Big\} \right) \\
 &\leq \sum_{\substack{(x_0,\dots, x_{n-1})\in L_{n-1} \\ e\in\ebd S_{\om,x_{n-1}}}}
\pr^\om\left( x_0\xlra[S,(x_0,\dots, x_{n-1})]{\om, p}x_{n-1}\right) p \, \pr^\om\bigg(x_{n-1}\xlra[S_{\om,x_{n-1}}]{\om, p} e^- \bigg)\\
 &= \sum_{(x_0,\dots, x_{n-1})\in L_{n-1}}
\pr^\om\left( x_0\xlra[S,(x_0,\dots, x_{n-1})]{\om, p}x_{n-1}\right)\phi_p^{\om,x_{n-1}}(S_{\om,x_{n-1}})
\leq \ex^\om\left(|L_{n-1}(p)|\right)c \, . 
\end{align*}
It follows by induction that $\ex^\om\big(|L_n(p)|\big)\leq c^n$. Therefore, 
\[\pr^\om(|\C_o|=\infty)=\pr^\om(T^\om(p) \textrm{ survives})=\lim_{n\to\infty}\pr^\om(|L_n(p)|\geq 1)\leq
\lim_{n\to\infty}\ex^\om|L_n(p)|=0\,, \]
hence $p\leq \pc$. 
 
Next we prove that $\pl\geq\pc$.
Let
\[\q(p):=\pp(\big\{\om: \phi^\om_p(S)\geq 1 \textrm{ for all } S\in\fss(\om)\}\big),\]
Note that $q(p)$ is non-decreasing in $p$, and $\q(p)>0$ for every $p>\pl$ by the definition of $\pl$. 

Fix $\om$ and let $H\in\fss(\om)$ be fixed.
We will use Lemma 1.4. of \cite{DCT}: 
\[\frac{d}{dp}\pr_p^\om\left(o\xlra{\om,p}H^c\right) 
\geq \left(1-\pr_p^\om\left( o\xlra{H,p} H^c \right)\right) \inf_{S: o\in S\subseteq H}\phi_p^H(S) 
\geq C(p) \inf_{S: o\in S\subseteq H}\phi_p^H(S)\,, \]
where $C(p)=(1-p)^{D} \leq 1-\pr^\om\left( o\xlra{\om,p} H^c \right)$ for every $\om$ and $H$, with
$D$ being the almost sure bound on the degree of the graph $G$. 
The probabilities above depend only on the structure of $\om$ in $K=H\cup \ovbd H$, hence
we can use the above inequality
to estimate the derivative of the probability $\pp\left(o\xlra{\om,p}B^\om(o,r)^c\right)$, as follows.
Consider the following sets of finite rooted graphs:
let $\mathcal{H}_r$ be the set of possible $(r+1)$-neighbourhoods of the
graphs with degree at most $D$, i.e.  
\begin{align*}
\mathcal{H}_r:=&\left\{(K,o): \mathrm{dist}_{K}(o,x)\leq r+1  \textrm{ and
}\mathrm{deg}_{K}(x)\leq D\,,\textrm{ for all }x\in V(K)\right\}, 
\end{align*}
and let
\begin{align*}
\mathcal{H}_r(p):=&\left\{(K,o)\in \mathcal{H}_r:  \phi_p^{K}(S)\geq 1\,,\textrm{ for all }S\in\fss(B_{K}(o,r))\right\}.
\end{align*}
Note that 
\begin{align*}
\sum_{K\in\mathcal{H}_r(p)}\pp\left(\left\{\om:
B_\om(o,r+1)=K\right\}\right)&=\\
\pp(\big\{\om: \phi^\om_p(S)\geq 1 \textrm{ for all } S\in\fss(B_{\om}(o,
r))\}\big)&\geq q(p),
\end{align*}
hence we have 
\begin{align*}
\frac{d}{dp}\pp\left(o\xlra{\om,p}B(o,r)^c\right) 
 &= \sum_{(K,o)\in\mathcal{H}_r} \pp\big(B_\om(o,r+1)=K\big)\frac{d}{dp}\pr_p\left(o\xlra{K,p}B_K(o,r)^c\right) \\
 &\geq \sum_{(K,o)\in\mathcal{H}_r(p)}\pp\big(B_\om(o,r+1)=K\big) C(p) \inf_{S: o\in S\subseteq B_{K}(o,r)}\phi_p^K(S)\\
 &\geq \q(p)\,C(p). 
\end{align*}
Integrate the above inequality on the interval $\left[\frac{p+\pl}{2},p\right]$.
Using the monotonicity of $q(p)$ and $C(p)$, we get
\[ \pp\left(o\xlra{\om,p}B(o,r)^c\right)
\geq \frac{p-\pl}{2}\, \q\left(\frac{p+\pl}{2}\right)C(p). \]
This gives a positive lower bound that is uniform in $r$. 
Thus $\pp\left(o\xlra{\om,p}\infty\right) >0$, and $p\geq \pc$.
\end{pfo}
\medskip

One advantage of the definition of $\pl$ for transitive graphs is that it enables one to check whether a certain 
$p$ is under $\pl$ using a finite witness. 
This characteristic makes the next definition natural. 

\begin{defi}\label{unifgood}
We say that a bounded degree unimodular random graph $G$ is \emph{uniformly good} if
for any $p<\pc$ there exists a positive integer $r(p)$ such that
$\pp_G(\{ \om: \exists S_\om\subseteq B_\om(o,r(p)), o\in S_\om \textrm{ s.t. }\phi_p^\om(S_\om)<1\})=1$ .
\end{defi}

This class of graphs includes unimodular quasi-transitive graphs (obvious) and unimodular random trees of uniform sub-exponential growth (see Definition \ref{subexp.def} and the proof of Proposition~\ref{subexp} in Subsection~\ref{ss.semi}). Furthermore, uniformly good unimodular graphs satisfy the following exponential decay of $\phi_p(B_\om(o,r))$ in $r$. 

\begin{lemma}\label{unifgoodexp}
Let $G$ be a bounded degree unimodular random graph. $G$ is uniformly good if and only if for all $p<\pc$ 
there are constants $c=c(p)<1$ and $R(p)$ such that if $r\geq R(p)$, then 
$\phi_p^\om(B)\leq c^r$ for almost every $\om$ and every finite
$B\supseteq B_\om(o,r)$.
\end{lemma}

For the proof of Lemma \ref{unifgoodexp} we use the same tree $T^\om$ as in the proof of Theorem \ref{pl=pc}. 
The uniformly good property implies a uniform linear lower bound in $r$ 
on the distance of the root from any vertex of $T^\om$ that corresponds to a boundary point of $B$ 
(namely the points of the set $\pi$ defined in the proof). 
This property and the boundedness of the size of the sets $S_{\om,x}$ allows us to prove the estimate of the lemma.
\medskip

\begin{pf}
If  the constants $c(p)$ and $R(p)$ exist, then the sets $S_\om:=B_\om(o,R(p))$ indicate that $G$ is uniformly good. 

To prove the other direction, assume that $G$ is uniformly good, and fix $p<\pc$. 
We can show as in the proof of Theorem \ref{pl=pc} that there exists a constant $c_0<1$ and a positive integer $r_0$
such that for almost every $\om$ and every $x\in\om$ there exists a finite connected set $S_{\om,x}\subseteq B_\om(x,r_0)$ containing $x$
that satisfies $\phi_p^{\om,x}(S_{\om,x})\leq c_0$. 
Fix an $\om$ and the sets $S_{\om,x}$ as above, a positive integer $r$ and a finite set $B\supseteq B_\om(o,r)$.
We define the trees $T^\om$ and $T^\om(p)$ as in the proof of Theorem \ref{pl=pc}.
On every directed path in $T^\om$ from $o$ to infinity there is a first vertex $(x_0,\dots ,x_k)$ such that $x_k\notin B$. 
Let $\pi$ be the set of these vertices, i.e.
\begin{align*}
\pi:=\{(x_0,\dots, x_k)\in T^\om: x_0,\dots, x_{k-1}\in B, x_k\notin B\}.
\end{align*}
Note that $\pi$ is a minimal set in $T^\om$ that separates $o$ from infinity, 
hence every non-backtracking infinite path from $o$ has exactly one vertex in $\pi$. 
An argument as in the first part of the proof of Theorem~\ref{pl=pc} shows that 
\begin{align*}
\ex^\om\left(|\pi\cap T^\om(p)|\right)&
=\sum_{(x_0, \dots ,x_k)\in \pi} 
\pr^\om\Big( x_0\xlra[B,(x_0, \dots ,x_k)]{\om, p} x_k\Big)  \\
 &\leq\sum_{(x_0, \dots ,x_k)\in \pi} \sum_{(x'_1,\dots, x'_k)}
\prod_{j=1}^k\pr^\om\Big( x_{j-1}\xlra[S_{\om,x_{j-1}}]{\om, p} x_{j}'\Big)p=:F(\pi,p), 
\end{align*}
where $(x'_1,\dots, x'_k)$ denotes a sequence of vertices in $\om$ such that $x_j'\in S_{\om, x_{j-1}}$ and $ x_j'\sim x_j$ for any $j=1,\dots, k$. 
First we bound $\phi_p^\om(B)$ in terms of $F(\pi,p)$ using the uniform bound on the size of the sets $S_{\om,x}$, 
then we prove a geometric bound on $F(\pi,p)$ using a linear bound in $r$ on the distance of $o$ and $\pi$ in $T^\om$. 
These two estimates will imply the statement of the lemma. 

Denote by $\bar{\pi}$ the set of the parents of the vertices in $\pi$, i.e. 
\begin{align*}
\bar{\pi}:=\{(x_0,\dots, x_k)\in T^\om: x_0,\dots, x_{k}\in B, \exists x_{k+1}\notin B, (x_0,\dots, x_{k+1})\in T^\om\}.
\end{align*}
If for some $e\in\ebd B$ the event $\big\{o\xlra[B]{\om,p}e^-\big\}$ occurs, then there is some $(x_0,\dots ,x_k)\in \bar{\pi}$ such that 
there is a good path through $x_0,\dots ,x_k$ in the percolation and a disjoint path from $x_k$ to $e^-$ in $S_{\om,x_k}$.
For any fixed $(x_0,\dots ,x_k)$ the number of edges in $\ebd B\cap \left(E(S_{\om,x_k})\cup \ebd S_{\om,x_k}\right)$ is bounded above by
$|E(S_{\om,x_k})\cup \ebd S_{\om,x_k}|\leq D^{r_0+1}$ where $D$ is the almost sure bound on the degree of the graph $G$. 
We have 
\begin{align}
\phi_p^\om(B)&=p\sum_{e\in\ebd B}\pr^\om\Big(o\xlra[B]{\om,p}e^-\Big) \nonumber \\
 &\leq \sum_{e\in\ebd B}\sum_{\substack{(x_0,\dots ,x_k)\in\bar{\pi}\\(x'_0,\dots ,x'_k)}}
\pr_p^\om\left(\Big\{x_0 \xlra[B,(x_0,\dots ,x_k)]{\om,p} x_k\Big\}\square\Big\{x_k\xlra[S_{\om,x_k}]{\om,p}e^-\Big\}\right) \nonumber\\
 &\leq \sum_{\substack{(x_0,\dots ,x_k)\in\bar{\pi}\\(x'_0,\dots ,x'_k)}} 
\bigg(\prod_{j=1}^k\pr^\om\Big( x_{j-1}\xlra[S_{\om,x_{j-1}}]{\om, p} x_{j}'\Big)p\bigg)
\sum_{e\in\ebd B\cap \left(E(S_{\om,x_k})\cup \ebd S_{\om,x_k}\right)}\pr_p^\om\Big(x_k\xlra[S_{\om,x_k}]{\om,p}e^-\Big) \nonumber \\
 &\leq \sum_{\substack{(x_0,\dots ,x_k)\in\bar{\pi}\\(x'_0,\dots ,x'_k)}}
\bigg(\prod_{j=1}^k\pr^\om\Big( x_{j-1}\xlra[S_{\om,x_{j-1}}]{\om, p} x_{j}'\Big)p\bigg)D^{r_0+1} 
=F(\bar{\pi},p)D^{r_0+1}. \label{phiF}
\end{align}
To estimate \eqref{phiF}, note that 
\begin{align*}
F(\pi,p)&= \sum_{\substack{(x_0,\dots ,x_k)\in\bar{\pi}\\(x'_0,\dots ,x'_k)}} 
\bigg(\prod_{j=1}^k\pr^\om\Big( x_{j-1}\xlra[S_{\om,x_{j-1}}]{\om, p} x_{j}'\Big)p\bigg) 
\sum_{\substack{x_{k+1}: (x_0,\dots ,x_{k+1})\in \pi \\ x'_{k+1}\in S_{\om,x_k}, x'_{k+1}\sim x_{k+1}}} 
\pr^\om\Big( x_{k}\xlra[S_{\om,x_{k}}]{\om, p} x_{k+1}'\Big)p\\
 &\geq \sum_{\substack{(x_0,\dots ,x_k)\in\bar{\pi}\\(x'_0,\dots ,x'_k)}} 
\bigg(\prod_{j=1}^k\pr^\om\Big( x_{j-1}\xlra[S_{\om,x_{j-1}}]{\om, p} x_{j}'\Big)p\bigg) 
p^{r_0+1}
= F(\bar{\pi},p)p^{r_0+1}
\end{align*}
 by the assumption that the graph is uniformly good. Combined this with \eqref{phiF} gives 
\begin{align}\label{F}
\phi_p^\om(B) 
 &\leq \frac{D^{r_0+1}}{p^{r_0+1}} F(\pi,p).
\end{align}

Now we show that $F(\pi,p)\leq c_0^{\frac{r}{r_0}}$, which combined with \eqref{F} proves the lemma.
Let $\pi_n:=\bigcup_{m\leq n}\left(\pi\cap L_m\right) \cup \left\{v\in L_n: v \textrm{ has a descendant in }\pi\right\}$, which is a minimal vertex set that separates the root from infinity.
Let $R:=\max\{n: L_n\cap\pi\neq \emptyset\}<\infty$, thus $\pi=\pi_R$. Note that each $\pi_n$ is the disjoint union of $\pi_{n+1}\setminus L_{n+1}\subseteq \pi$ and $\pi_n\setminus \pi_{n+1}\subseteq L_n$. 
We estimate $F(\pi,p)$ by summing over a larger set: the union of
$\pi_R\setminus L_R$ and $\{(x_0, \dots ,x_{R}): (x_0, \dots ,x_{R-1})\in
\pi_{R-1}\setminus \pi_{R}, x_R\in \ovbd S_{\om, x_{R-1}}\}\supseteq \pi_R\cap
L_R$. That is, using the bound 
\begin{align*}
\sum_{e\in \ebd S_{\om, x_{R-1}}} \pr^\om\Big(
x_{R-1}\xlra[S_{\om,x_{R-1}}]{\om, p} e^- \Big)p=\phi_p^{\om,
  x_{R-1}}(S_{\om, x_{R-1}})\leq 1
\end{align*} 
for the second term in the following estimation, we have 
\begin{align*}
F(\pi, p) & \leq\sum_{(x_0, \dots ,x_k)\in \pi_R\setminus L_R} \sum_{(x'_1,\dots, x'_k)}
\prod_{j=1}^k\pr^\om\Big( x_{j-1}\xlra[S_{\om,x_{j-1}}]{\om, p} x_{j}'\Big)p \\
& +\sum_{\substack{(x_0, \dots ,x_{R-1})\in \pi_{R-1}\setminus \pi_{R} \\ (x'_1,\dots, x'_{R-1})}} 
\bigg(\prod_{j=1}^{R-1}\pr^\om\Big( x_{j-1}\xlra[S_{\om,x_{j-1}}]{\om, p} x_{j}'\Big)p\bigg)
\sum_{e\in \ebd S_{\om, x_{R-1}}} \pr^\om\Big( x_{R-1}\xlra[S_{\om, x_{R-1}}]{\om, p} e^- \Big)p \\
& \leq\sum_{(x_0, \dots ,x_k)\in \pi_{R-1}} \sum_{(x'_1,\dots, x'_k)}
\prod_{j=1}^k\pr^\om\Big( x_{j-1}\xlra[S_{\om,x_{j-1}}]{\om, p} x_{j}'\Big)p = F(\pi_{R-1},p).
\end{align*}
A similar argument shows that $F(\pi,p)\leq F(\pi_{n},p)$ for any $n\leq R$. 
If $(x_0,\dots ,x_k)\in\pi$, then $\dist_\om(o, x_k)\geq r$, hence 
the distance between $o$ and $\pi$ in $T^\om$ is at least $\frac{r}{r_0}$, thus $\pi_{n}=L_{n}$ for any $n\leq \frac{r}{r_0}$. 
If we apply the above argument for $F(\pi_{n},p)$ whith $n\leq \frac{r}{r_0}$, then the first term disappear, and the inequality $\phi_p^{\om, x_{n-1}}(S_{\om, x_{n-1}})\leq c_0$  gives  
\begin{align*}
F(\pi,p)\leq F(\pi_{\frac{r}{r_0}},p)\leq F(\pi_{\frac{r}{r_0}-1},p)c_0 \leq \dots \leq c_0^{\frac{r}{r_0}}\,.  
\end{align*}
This combined with \eqref{F} proves the lemma. 
\end{pf}

\begin{cor}\label{pcpe}
If $G$ is a uniformly good unimodular graph, then $\pc\leq\pe$. 
\end{cor}

\begin{pf}
Let $p<\pc$, and let $c$ and $R(p)$ be as in Lemma \ref{unifgoodexp}. We have $\ee\left(\phi_p^\om\left(B_\om(o,R(p))\right)\right)\leq c^R<1$, thus $p\leq \pe$.
\end{pf}
\medskip

We will see in Remark \ref{pc>pe} that, without the assumption of uniform goodness, 
the inequality $\pc\leq \pe$ does not necessarily hold.  
Also, we will show in Example \ref{pc<pe} that there are uniformly good graphs with $\pc<\pe$.  

\begin{defi}\label{subexp.def}
A random rooted graph $(G,o)$ has \it{uniform sub-exponential volume growth} 
if for any $c<1$ and $\eps>0$ there is an $R$ such that 
$\pp\left(\om:  |B_\om(o,r)|c^r<\eps\right)=1$ for any $r>R$. 
\end{defi}

\begin{cor}\label{pc=pt} 
If $G$ is a uniformly good unimodular graph with uniform sub-exponential volume growth, then $\pc=\pt=\ptk$.
\end{cor}

\begin{pf}
Let $p<\pc=\pl$ and let $c$ and $R(p)$ be as in Lemma \ref{unifgoodexp}. 
Denote by $D$ the maximum degree of $G$. 
Let $R>R(p)$ such that $\pp\left(\{\om:  |B_\om(o,r)|c^{r/2}<1 \}\right)=1$ for any $r>R$ and 
let $\om$ satisfy this event for all $r>R$ simultaneously. 
%\in\bigcap_{r>R}\left\{\om': |B_{\om'}(o,r)|c^{r/2}<1 \right\}$. 
Then we have 
\begin{align*}
\ex_p^\om\left( |\C_o| \right) &
= \sum_{n=1}^{\infty} \pr^\om_p\left( |\C_o|\geq n \right) 
= \sum_{r=1}^{\infty}\sum_{n=|B_\om(o,r)|+1}^{|B_\om(o,r+1)|} 
\pr^\om_p\left( |\C_o|\geq n \right) \\
 &\leq \sum_{r=1}^{\infty}\sum_{n=|B_\om(o,r)|+1}^{|B_\om(o,r+1)|} 
\pr^\om_p\left( o\xlra{p,\om} B_\om(o,r)^c \right) \\
 &\leq \sum_{r=1}^{\infty} |B_\om(o,r+1)| \min\{ \phi^\om_p\left( B_\om(o,r) \right) , 1 \}\\
 &\leq  \sum_{r=2}^{R+1} |B_\om(o,r)|+ \sum_{r=R+1}^{\infty} |B_\om(o,r+1)|c^r \\
 &\leq \sum_{r=2}^{R+1} D^r+ \sum_{r=R+1}^{\infty} c^{r/2} <\infty
\end{align*}
This gives a uniform upper bound on $\ex_p^\om\left( |\C_o| \right)$ 
thus $\ee_p\left( |\C_o| \right)<\infty$. 
It follows that $p\leq \ptk$, hence $\ptk \geq \pc$. 
The other direction follows from the definition of $\ptk$. 
\end{pf}
\medskip

Subexponential volume growth also appears in Theorem~\ref{subexp}, Example~\ref{Zd} and Corollary~\ref{biLip}.

\subsection{Counterexamples}\label{ss.pc.countex}

We show in Examples~\ref{subexpcanopy} and~\ref{ptk} that there are unimodular random graphs of
uniform subexponential (in fact, quadratic) volume growth, but $\pt<\pc$ and $\ptk <\pt$. Both
constructions will use Bernoulli percolation on $\Z^2$ as an ingredient;
moreover, although we define the graph in the second  example as a vertex
replacement of $\Z^2$, it could be defined even as an invariant random
subgraph of $\Z^2$. We further give examples of graphs with $\pe<\pc$ and
$\pe>\pc$; see Examples~\ref{pc>pe} and~\ref{pc<pe}, respectively. 
First we need a lemma that will be useful in our examples. 

\begin{lemma}\label{crossing}
For any $\eps>0$ there is a probability $p_1<1$ such that for $n$ large enough,
the vertices $(0,-n)$, $(0,n)$, $(-n,0)$, $(n,0)$ are in the same cluster in Bernoulli($p_1$) percolation on $Q_n$
with probability at least $1-\eps$.
\end{lemma}

\begin{pf} 
The occurrence of the events in the following two claims implies the
occurrence of the event in the statement of the lemma, hence we will be done
by a union bound. 
\medskip

\noindent {\it Claim 1:} For any $p>{1}/{2}$ and $n>n_0(p,\eps)$ large enough, in Bernoulli($p$) percolation on $Q_n$, with probability at least $1-\eps/2$, there is a \emph{giant cluster} with the following properties: it joins all the sides of $Q_n$, while every other cluster in $Q_n$ has diameter at most ${n}/{5}$. This was proved in \cite[Proposition 2.1]{AP}.
\medskip

\noindent {\it Claim 2:} There exists $p_1<1$ such that for all $n$ and all $p>p_1$,
$$
\pr_p\left(\mathrm{diam}(\C_{(0,n)}) \geq n\right) \geq 1-\eps/8\,.
$$
Similarly for $(0,-n)$, $(-n,0)$, and $(n,0)$, instead of $(0,n)$. The proof follows from a standard Peierls contour argument, thus we leave it to the reader.
\end{pf}
\medskip

We will use the following unimodular random graph, the canopy tree, in several of our examples. It is the local weak limit of large balls in the 3-regular tree: 

\begin{defi}[Busemann functions and canopy tree]\label{canopytree}
Let $\T$ be the 3-regular infinite tree with a root $o$, a distinguished end $\xi$, and 
a \emph{Busemann function} (see \cite{W}) $\mathfrak{h}:\T\to \Z$ that gives the levels w.r.t.~to $\xi$. 
More precisely, to define $\mathfrak{h}$, for any vertex $x$, let $(\xi,x)$ be the unique infinite simple path from $x$ which is in the equivalence class $\xi$. Denote by $x\wedge o$ the unique vertex in $\T$ such that $(\xi,x\wedge o)=(\xi,x)\cap (\xi,o)$, 
i.e., the first vertex where $(\xi,x)$ and $(\xi,o)$ coalesce. Finally,
let $\mathfrak{h}(x):=\dist(o,x\wedge o)-\dist(x,x\wedge o)$. 

Let $\Lambda\subset\T$ be the subgraph spanned by the vertices $x$ with 
$\mathfrak{h}(x)\geq 0$. This tree $\Lambda$ is called the \emph{canopy tree}. Denote by $L(n):=\{x\in V(\T):\mathfrak{h}(x)=n\}$ the $n^\textrm{th}$ vertex level and by $L_E(n):=\{e\in E(\T):e^-\in L(n), e^+\in L(n+1)\}$ the $n^\textrm{th}$ edge level of $\T$, or, for $n\geq 0$, of $\Lambda$. If we choose the root $o$ of $\Lambda$ such that $\pr(o\in L(n))=2^{-n-1}$, we get a unimodular random graph.
\end{defi}

\begin{exa}\label{subexpcanopy}
There is a unimodular graph with uniform polynomial volume growth and $\pt<\pc$. In particular, the exponential decay of two-point connection probabilities fails for $p\in (\pt,\pc)$ on this graph.
\end{exa}

\begin{pf}
We define the graph $G$ as an edge replacement (see \cite{AL07}, Example 9.8) of the canopy tree: each $e\in L_E(n)$ is replaced by $(Q_{2^n}(e),(0,-{2^n}),(0,{2^n}))$, where $Q_{2^n}(e)$ is isomorphic to $Q_{2^n}$. It is easy to see that the volume of $B_G(o,r)$, for any root $o$ and radius $r$, is at most $C r^2$, for some absolute constant $C < \infty$. Indeed, if the root is in $Q_{2^n}(e)$, then $B_G(o,r)$ intersects the cubes $Q_{2^l}(e')$ with $e'\in(\xi,e)$ only if $l\leq \log_2 r$ or $l=n$. Furthermore, each such $Q_{2^l}(e')$ has more vertices than the sum of the number of vertices of $Q_{2^k}(e'')$ with $e'\in (\xi,e'')$, which are the further cubes that may intersect $B_G(o,r)$. It follows that 
$|B_G(o,r)|\leq \min\left\{r^2, \sum_{l=n}^{\log_2 r} 2^{2l+3}\right\}\leq Cr^2$. 

We will now show that $\pt(G) < \pc(G)=1$. Consider Bernoulli($p$) percolation $\omega$ on $G$ and, as a deterministic function of it, define the following percolation $\lambda$ on $\Lambda$: an edge $e\in L_E(n)$ is open in $\lambda$ if and only if the vertices $(0,-n)$ and $(0,n)\in Q_n(e)$ are connected by an open path in $\omega$. Clearly, there exists an infinite cluster in $\omega$ if and only if there is an infinite cluster in $\lambda$. The law of $\lambda$ is stochastically dominated by a Bernoulli($1-(1-p)^3$) percolation on $\Lambda$, because if $e\in L_E(n)$ is open, then at least one of the edges in $Q_n(e)$ adjacent to $(0,n)$ is open. The tree $\Lambda$ has one end, hence, for any $p<1$, 
$$
\pr_p^G(\exists \text{ an infinite cluster})\leq \pr_{1-(1-p)^3}^\Lambda(\exists \text{ an infinite cluster})=0\,.
$$ 
That is, $\pc(G)=1$.

An easy first moment computation (that we omit) shows that $\pt(\Lambda)={1}/{\sqrt{2}}$. Now let $0<\eps <1-{1}/{\sqrt{2}}$. 
It follows from Lemma~\ref{crossing} that there exists $p_1<1$ and some large $N$ such that 
$\pr_{p_1}(e\in\lambda)\geq 1-\eps$ for all $e\in L_E(n)$ with $n\geq N$. 
Thus, for $o\in L(N)$, the cluster $\C_o$ in $\lambda$, restricted to the levels $n\geq N$, stochastically dominates Bernoulli($1-\eps$) percolation on $\Lambda$. The latter has infinite expected size, hence the expected size of the cluster in $\omega$ of $(0,-N)\in Q_N(e)$ for $e\in L_E(N)$ is also infinite. That is, $\pt(G)\leq p_1<1$. 
\end{pf}

\begin{exa}\label{pc>pe}
The canopy tree $\Lambda$ (see Definition \ref{canopytree}) satisfies
$\pe=\frac{1}{\sqrt{2}}$, thus this is an example of a not uniformly good
unimodular graph with $\pc>\pe$. 
\end{exa}

\begin{pf}
It is easy to check that $\ee(\phi_p(B(o,r)))$ equals $2p(\sqrt{2}p)^r$ if $r$ is even, and equals
$3(\sqrt{2}p)^{r+1}/2$ if $r$ is odd. This sequence converges to 0 for $p<1/\sqrt{2}$, while remains above 1 for $p>1/\sqrt{2}$, which implies the claim.
\end{pf}

\begin{exa}\label{ptk}
There is a unimodular graph with polynomial volume growth and $\ptk <\pt$.
\end{exa}

\begin{pf}
Let $X$ be a positive integer valued random variable such that $\pr(X=k)=ck^{-5/2}$ for all $k\ge 1$. 
Then $\ex X<\infty$ and $\ex (X^2)=\infty$. 
We define the graph $G$ as a vertex replacement (see Subsection~\ref{ss.opera}) of $\Z^2$ with respect to the following labels as follow. Let $\{X_n, X'_n: n\in \Z\}$ be iid copies of $X$, and for each vertex $(m,n)\in\Z^2$, let 
$G_{(m,n)}$ be isomorphic to the subgraph of $\Z^2$ spanned by the vertices in $[0,2X_m]\times [0,2X'_n]$, and for the edges going from $(m,n)$ to North, East, South, and West, let the image of $\varphi_{(m,n)}$ be the corresponding midpoint of the box $G_{(m,n)}$. We can also think of the resulting graph as an invariant random subgraph of $\Z^2$.

Denote by $Y$ and $Y'$ half the length of the sides of the box of $o$ in $G$, i.e., 
the law of $X_0$ and $X'_0$ biased by $X_0X'_0$. Then 
\[ \pr(Y=k, Y'=l)=\frac{kl}{(\ex X)^2}\pr(X=k, X'=l), \]
hence $Y$ and $Y'$ are independent with distribution $\pr(Y=k)=\frac{ck^{-3/2}}{\ex X}$. 

First we show that $\ptk=\frac{1}{2}$. 
$G$ is a subgraph of $\Z^2$, hence $\ptk(G)\geq \frac{1}{2}$.  
Fix $p>\frac{1}{2}$ and let $\eps >0$. 
Denote by $M(Q_n)$
the largest cluster in percolation with parameter $p$ in the box $Q_n$, and let  
\begin{align*}
\A(Q_n):=\big\{|M(Q_n)|\geq (1-\eps)\theta(p)|Q_n|,\quad \mathrm{diam}(C)< \nu\log n
\ \forall 
\textrm{ open cluster } C\neq M(Q_n) \big\},
\end{align*} 
where $\theta(p)=\pr_p(|\C_o(\Z^2)|=\infty)$, and $\nu$ is chosen as follows:
by \cite[Theorem 7.61]{Gr}, there is an $N=N(p)$ and $\nu=\nu(p)$ such that, for any $n\geq N$,
\begin{align*}
\pr_p(\A(Q_n))>1-\eps\,.
\end{align*}
Let $Z:=\min\{Y, Y'\}$, and consider the event $\mathcal{D}(G_{0,0}):=\left\{\dist(o, \ivbd G_{0,0})\geq\nu\log Z\right\}$. 
If $ Z$ is large enough, then $\pr\big(\mathcal{D}(G_{0,0})\,\big| \,  Z \big)\geq
1-\eps$, since $o$ is uniform in $G_{0,0}$. Assuming that $\mathcal{D}(G_{0,0})$ occurs, choose a box
$Q_{ Z}\subseteq G_{0,0}$ that
contains $o$ such that $\dist(o,\ivbd Q_{ Z})\geq \nu\log  Z$. Consider percolation on $\Z^2\supset
Q_{ Z}$. If $o$ is in the unique infinite cluster of this percolation
on $\Z^2$, then the
diameter of $\C_o(Q_{ Z})$ is at least $\nu\log Z$, hence 
\begin{align*}
\pr_p\Big(o\in M(Q_{ Z}), \A(Q_{ Z})\,\Big|\, Z=n, \mathcal{D}(G_{0,0})\Big)>\theta(p)-\eps
\end{align*}
for $n$ large enough. It follows that there is an $N'$ such that 
\begin{align*}
\ee_p\left(|\C_o|\right)& \geq \sum_{n=N'}^{\infty}\pr_p\Big(o\in M(Q_{ Z}),
\A(Q_{ Z}), \mathcal{D}(G_{0,0})\Big| Z=n \Big) \, \pr( Z=n)\,(1-\eps)\theta(p)n^2 \\ 
& \geq \sum_{n=N'}^{\infty}(\theta(p)-\eps)(1-\eps) \, \pr( Z=n) \, (1-\eps)\theta(p)n^2  =\infty\,, 
\end{align*}
as desired.

To show that $\pt>\frac{1}{2}$ let $e$ be an edge in $\Z^2$, and let $G_{e^-}$ and $G_{e^+}$ be the subgraphs of $G$ that correspond to the endpoints of the edge. Let $x:=\varphi_{e^-}(e)$ and $y:=\varphi_{e^+}(e)$, i.e. let $\{x,y\}$ be the edge in $G$ that joins $G_{e^-}$ and $G_{e^+}$. 
If there is an open path in $G(p)$ through the edge $\{x,y\}$, that joins two vertices in $G_{e^-}\setminus \{x\}$ and in $G_{e^+}\setminus \{y\}$, then the event $J(\{x,y\}):=\left\{\exists e'\in E(G_{e^-}): e'\sim x, e' \textrm{ open}\right\}\cap \left\{\exists e'\in E(G_{e^+}): e'\sim y, e' \textrm{ open}\right\} \cap \left\{\{x,y\} \textrm{ open}\right\}$ occurs. 
For a fixed configuration of $G$ the events $J(\{\varphi_{e^-}(e),\varphi_{e^+}(e)\})$ are independent for different edges, and $\pr_p(J(\{\varphi_{e^-}(e),\varphi_{e^+}(e)\}))=p(1-(1-p)^3)^2$. 
This probability is strictly increasing in $p$ and there is a $p_0>\frac{1}{2}$ such that 
$p(1-(1-p)^3)^2>\frac{1}{2}$ if{f} $p>p_0$. 
We consider a random subset $H=H\left(G(p)\right)\subseteq E(\Z^2)$ obtained from the percolation $G(p)$: 
let $e\in H$ if and only if the event $J(\{\varphi_{e^-}(e),\varphi_{e^+}(e)\})$ occurs in $G(p)$. 
The law of $H$ is the same as the law of Bernoulli($p(1-(1-p)^3)^2$) bond percolation. 
We want to estimate the expected size of $\C_o(G)$ conditioned on the size of $G_{0,0}$. 
If $\C_o(G)$ intersects a box $G_v$, then the connected component of $o$ in $H$ contains $v$. 
Therefore 
\begin{align*}
\ee_p\Big( |\C_o| \,\Big| \,Y, Y' \Big)&
\leq  \ee_p\left(\left.\sum_{v\in \Z^2: v\in \C_o(H)} |G_v|\right|Y,Y'\right)  \\
 &\leq \ex(|\C_o(H)|)\max\left\{Y^2, (Y')^2, (\ex X)^2\right\}, 
\end{align*}
which is finite if $p<p_0$. 
It follows that for almost every configuration of $(G,o)$ the expected size $\ex^{G}_p(\C_o)$ is finite if $p<p_0$, 
hence $\pt\geq p_0$.  
\end{pf}

\begin{exa}\label{pc<pe}
There is a quasi-transitive graph with $\pe>\pc$.
\end{exa}
\begin{pf}
Let $H_{k,l}$ be the following finite directed multigraph: the vertex set is $\{x_0,x_1,\dots, x_k\}$,
 and we have $l$ loops at $x_0$, then one edge from
$x_0$ to each $x_j$, $j=1,\dots,k$, and one from each $x_j$ back to $x_0$. 
Let $T_{k,l}$ be the directed cover of $H_{k,l}$ based at $x_0$. Consider two copies of
$T_{k,l}$ and connect the roots of them by an edge to get the infinite quasi-transitive graph $G_{k,l}$,
which has vertices of degree 2 and $k+l+1$.  
One can easily compute that to get a unimodular random graph one has to
choose the root according to  $\pp(\mathrm{deg}\,o=2)=1-\pp(\mathrm{deg}\,o=k+l+1)=\frac{k}{k+2}$.
Hence $\ee(\mathrm{deg}\,o)=\frac{4k+2l+2}{k+2}$.
The equality $\ee\left(\phi_p^\om(B_\om(o,0))\right)=p\,\ee(\mathrm{deg}\,o)$
implies that $\pe\geq\left(\ee(\mathrm{deg}\,o)\right)^{-1}=\frac{k+2}{4k+2l+2}$. 
On the other hand, the critical probability of a directed cover of a finite
graph is
$\pc(T_{k,l})=\left(\mathrm{br}(T_{k,l})\right)^{-1}=\left(\mathrm{growth}(T_{k,l})\right)^{-1}=\left(\lambda_*(H_{k,l})\right)^{-1}$,
where $\lambda_*(H)$ is the largest positive eigenvalue of the directed
adjacency matrix of $H_{k,l}$; see \cite{LP}, Section 3.3 and \cite{Lyo90}. One can thus compute that
$\pc(G_{k,l})=\pc(T_{k,l})=\frac{2}{l+\sqrt{l^2+4k}}$. If we set, e.g., $k=3,
l=5$, then we have 
$\pc(G_{3,5})=\frac{2}{5+\sqrt{37}}<\frac{5}{24}=\left(\ee_{G_{3,5}}(\mathrm{deg}\,o)\right)^{-1}\leq \pe(G_{3,5})$.
\end{pf}

\section{Locality of the critical probability}\label{s.loc}

In this section we examine the question of Schramm's locality conjecture: does $\pc(G_n)$ converge to $\pc(G)$ if $G_n\to G$ in the local weak sense? 
The original question in \cite{BNP} was phrased for sequences of transitive graphs that converge to a transitive graph in the local sense and satisfy $\sup \pc(G_n)<1$. 
First we provide some simple examples of unimodular graphs where the conjecture holds. 
In Example~\ref{perccluster}, we note that if $G_n$ and $G$ are infinite clusters of an independent percolation with appropriate parameters, then the convergence holds. In Example~\ref{ugwt}, we discuss unimodular Galton--Watson trees, and give sufficient and necessary conditions on the offspring distribution to satisfy locality of $\pc$. 
Then we investigate the inequality $\liminf\pc(G_n)\geq \pc(G)$, which is
known for transitive graphs; see \cite{DCT} for a simple proof, or the first paragraph of Subsection~\ref{ss.semi}.
%\gabor{inkabb ide kellene a bizonyitas, vagy pointer a Subsectionre.} 
In Proposition
\ref{liminfpe} we show by a similar argument that the critical probability
$\pe$ satisfies this inequality for unimodular random graphs. We prove Propositions~\ref{liminf} and~\ref{subexp} that state that 
under certain restrictions on the graphs $G$ and $G_n$ the convergence $\lim
\pc(G_n)=\pc(G)$ is true for unimodular random graphs.
Examples~\ref{lamplighter} and~\ref{Zd} provide graph sequences with $\lim\pc(G_n)<\pc(G)$.
These indicate that unimodular graphs do not satisfy Schramm's conjecture in general and show that the conditions in Proposition~\ref{liminf} and~\ref{subexp} are necessary. We show in Example~\ref{box} a sequence with $\pc(G)<\lim\pc(G_n)<1$. 
In this example $G$ and each $G_n$ satisfy the conditions of Corollaries~\ref{pcpe} and~\ref{pc=pt}, thus $\pc=\pt=\ptk$ and also $\pe(G)<\lim\pe(G_n)<1$. This shows that none of the generalisations of the critical probabilities satisfies the extension of Schramm's conjecture for unimodular graphs in general.

\subsection{Basic examples}\label{ss.basic}

We present now two natural classes of unimodular random graphs that satisfy Schramm's conjecture. The first example is very easy; the proof is left as an exercise. 

\begin{exa}\label{perccluster}
Let $\Gamma$ be a transitive unimodular graph and let $p_n\to p\in (\pc(\Gamma),1]$. 
Let $G_n$ (resp. $G$) be the connected component of the root in the Bernoulli($p_n$) (resp. $p$) percolation on $\Gamma$ conditioned to be infinite. 
Then $\pc(G_n)\to \pc(G)<1$. 
\end{exa}

Our second class of examples, unimodular Galton--Watson trees, is less trivial. 
Let $X$ be a non-negative integer valued random variable, the \emph{offspring distribution} of the tree, and 
let $UGW(X)$ be the unimodular Galton--Watson tree measure on rooted trees: 
the probability that the root $o$ has $k$ children is
\begin{equation}\label{UGW}
\pr^{UGW(X)}(\deg o=k)=\frac{\pr(X=k-1)}{k\,\ex(\frac{1}{X+1})}
\end{equation}
for $k\geq 1$, while the number of children of each descendant is according to $X$, independently of the other vertices. 
This measure is unimodular (see \cite{AL07}, Example 1.1), and if
$\ex X>1$, then $\pr(|UGW(X)|=\infty)>0$, thus we can consider the measure
$UGW_\infty(X)$ which is $UGW(X)$ conditioned on the event
$\{|UGW(X)|=\infty\}$. The measure $UGW_\infty$ is also unimodular, being an ergodic component of a unimodular measure.

\begin{exa}\label{ugwt}
Let $UGW_\infty(X)$ be the unimodular Galton--Watson tree with offspring distribution
$X$, conditioned to be
infinite. If $X_n$ and $X$ are non-negative integer valued random variables s.t.~the $X_n$ satisfy $\ex X_n>1$, while $X$ satisfies $\ex X>1$ or $\pr(X=1)=1$, then
\begin{itemize}
\item[{\bf (1)}]
$UGW_\infty(X_n)\to UGW_\infty(X)$ in the local weak sense if{f} $X_n\to X$ in distribution;
\item[{\bf (2)}] $\pc(UGW_\infty(X_n))\to \pc(UGW_\infty(X))$ if{f} $\ex X_n\to \ex X$. 
\end{itemize}
\end{exa}

Before the proof, note that this example shows that $\pc$ is a continuous function of $UGW_\infty(X)$ when the trees have a uniform bound on their degrees (by the Dominated Convergence Theorem), but not necessarily otherwise: if $X_n\to X$ in distribution, with $\ex X_n>1$ and $\ex X>1$, but $\ex X_n\nrightarrow \ex X$, then the critical probabilities $\pc(UGW_\infty(X_n))$ do not converge to $\pc(UGW_\infty(X))$. Nevertheless, Fatou's lemma implies that the inequality $\limsup \pc(UGW_\infty(X_n))\leq \pc(UGW_\infty(X))$ does hold without any assumptions. That is, if the trees do not satisfy the locality of $\pc$, then they also fail to satisfy the lower semicontinuity discussed in the next subsection, proved to hold in many cases, including transitive graphs. This suggests that a uniform bound on the degrees is a natural condition when we investigate the locality of $p_c$ for unimodular graphs. 
\medskip

\begin{pf} The critical probability $\pc(UGW_\infty(X))$ equals $\frac{1}{\ex X}$ (see \cite{LP}, Proposition 5.9), therefore $\pc(UGW_\infty(X_n))\to \pc(UGW_\infty(X))$ iff $\ex X_n\to \ex X$. This shows part (2).

For part (1), for any nonnegative integer random variable $X$, let $p_k(X):=\pr(X=k)$, 
let $f_X(t):=\sum_{k=0}^\infty p_k(X)t^k$ be the probability generating function of $X$, 
and let $q=q(X):=\pr(|GW(X)|<\infty)$, which is the smallest non-negative number that satisfies $f_X(q)=q$.

Assume that $X_n\to X$ in distribution, first with $\ex(X)>1$. From
$X_n\to X$ it follows easily that $UGW(X_n)\to UGW(X)$, while, from the uniform
convergence of the convex functions $f_{X_n}$ to the strictly convex function $f_X$ on $[0,1]$,
we also get $q_n=q(X_n)\to q(X) <1$. Thus $UGW_\infty(X_n) \to UGW_\infty(X)$.

Now assume that $\pr(X=1)=1$ and $\pr(X_n=1)\to 1$ with $\ex(X_n)>1$. Using~Bayes' rule and~(\ref{UGW}),
\begin{align}
\pr^{UGW_\infty(X_n)}(\deg o=2)
&=\frac{\pr^{UGW(X_n)}\left(|UGW(X_n)|=\infty\left|\deg o=2\right.\right)\pr^{UGW(X_n)}\left(\deg o=2\right)}{\pr\left(|UGW(X_n)|=\infty\right)}\nonumber \\
&=\frac{1-q_n^2}{\pr(|UGW(X_n)|=\infty)}\,\frac{\pr(X_n=1)}{2\,\ex(\frac{1}{X_n+1})}\nonumber\\
&=\frac{(1-q_n^2) \pr(X_n=1)}{2\sum_{j=1}^\infty \pr(X_n=j-1)(1-q_n^j) / j}\,. \label{qn}
\end{align}
We claim that $\pr^{UGW_\infty(X_n)}(\deg o=2)\to 1$. If $q_n$ converges to some $q_\infty < 1$, then plugging $\pr(X_n=1)\to 1$ into~\eqref{qn} yields the claim immediately. If $q_n\to 1$, then, simplifying the numerator and the denominator of~\eqref{qn} by $1-q_n$, it becomes
\begin{equation}\label{qn1}
\frac{(1+q_n) \, \pr(X_n=1)}{ 2 \sum_{j=1}^\infty \pr(X_n=j-1) (1+q_n+\dots+q_n^{j-1})/j} 
\geq \frac{(1+q_n) \pr(X_n=1)}{2} \to 1\,.
\end{equation}
Finally, if $q_n$ does not converge, we  can still apply one of these two arguments to any convergent subsequence, and obtain the claim. Therefore, in the local weak limit, the root has degree 2 almost surely. By unimodularity, every vertex has degree 2 almost surely (see \cite{AL07}, Lemma 2.3), hence this limit must be $\Z$. This is also $UGW_\infty(X)$, thus we have $UGW_\infty(X_n) \to UGW_\infty(X)$.

For the other direction of part (1), suppose that there are $X_n$ and
$X$ such that $UGW_\infty(X_n) \to UGW_\infty(X)$, but $X_n\nrightarrow
X$. The set $\{X_n\}$ of probability distributions must be tight: otherwise, 
a uniform random neighbour of $o$ in $UGW_\infty(X_n)$, 
whose offspring distribution stochastically dominates $X_n$ because of the conditioning on  $\big\{ | UGW(X_n) |=\infty\big\}$,
would have arbitrarily large degrees with a uniform positive probability, 
and thus $UGW_\infty(X_n)$ could not
converge to the locally finite graph $UGW_\infty(X)$. 
It follows from this tightness that there is a subsequence $\{X_{k(n)}\}$ that
converges in distribution to a random variable $Y\not=X$.

First we show that $\ex Y\geq 1$. Suppose $\ex Y<1$, then $\lim q_n =q(Y)=1$,
hence 
\begin{align*}
\pr^{UGW_\infty(X_n)}(\deg o=k)=\frac{\pr(X_n=k-1)(1+\dots
  +q_n^{k-1})}{k\sum_{j=1}^\infty \pr(X_n=j-1)(1+\dots +q_n^{j-1})/j}\to \pr(Y=k-1).
\end{align*}
It follows, that the expected degree of the root in the limit graph is $\ex
Y+1<2$. The local weak limit of the graphs $UGW_\infty(X_n)$ is almost surely
infinite, hence the expected degree of the root is at least 2 (see \cite{AL07},
Theorem 6.1), a contradiction. 

If we have $\pr(Y=1)=1$, then the first direction of part (1) implies that $UGW_\infty(X_{k(n)}) \to UGW_\infty(Y) = \Z$. But we also have  $UGW_\infty(X_{k(n)}) \to UGW_\infty(X)$, and it is obvious that $UGW_\infty(X)=\Z$ implies that $\pr(X=1)=1$. That is, $X_n$ would in fact converge in distribution to $X$, a contradiction.

If $\ex Y=1$, but $\pr(Y=1)\not=1$, then the generating function $f_Y(t)$ is
strictly convex, hence $q(X_{k(n)})\to q(Y)=1$. A computation similar to~\eqref{qn} and~\eqref{qn1} 
gives that the degree distribution of $o$ in $UGW_\infty(X_{k(n)})$ converges to that of $Y+1$. 
This must be the degree distribution of $o$ in the local limit $UGW_\infty(X)$. 
Since $\pr(Y+1=2)\not=1$, we must be in the case $\ex X>1$. 
However, then we would have $p_c(UGW_\infty(X))=1/\ex X<1$, while
$\ex(\deg o)=\ex(Y+1)=2$ implies that $UGW_\infty(X)$ is a tree with at most two ends (see \cite{AL07}, Theorem 6.2) hence $p_c=1$,
again a contradiction.

The final case is that $\ex Y>1$, for which we can again use the first direction of part (1), saying that $UGW_\infty(X_{k(n)})\to UGW_\infty(Y)$. If we prove that the distribution of $UGW_\infty(X)$ determines $X$, then we must have $X=Y$, and we are done, as before.
 
This invertibility follows from the construction in \cite{LP}, Theorem 5.28, as follows. Let $T^*:=GW(X^*)$, where the probability
generating function of the positive integer valued random variable $X^*$ is $f^*(t):=\frac{f_X(q+(1-q)t)}{1-q}$, and let
$\bar{T}:=GW(\bar{X})$, where
$\bar{f}(t)=f_{\bar{X}}(t):=\frac{f(qt)}{q}$, and hence $\bar{T}$ is almost
surely finite. The law of $GW(X)$ conditioned to be infinite 
equals the law of the tree $T$ constructed as follows: consider the rooted tree
$T^*$, and attach to each vertex of $T^*$ an
appropriate number of independent copies of $\bar{T}$. We get the
law of $UGW_\infty(X)$ if we attach to the root an appropriate random number
of independent copies of $T$ and $\bar{T}$. It follows that the law of
$UGW_\infty(X)$ determines $(f^*,\bar{f})$. We get the function $f$ from $(f^*,\bar{f})$ by the transform
$f(s)=q\bar{f}\left(\frac{s}{q}\right)$, if $0\leq s\leq q$ and
$f(s)=(1-q)f^*\left(\frac{s-q}{1-q}\right)$, if $q\leq s\leq 1$. There is a
unique $q$ for which the resulting $f(s)$ has the same second derivative from the left and from the right at $s=q$.
Since $f(s)$ has to be analytic, we see that $(f^*,\bar{f})$ uniquely determines $f$ and hence $X$. 
\end{pf}

\subsection{Semicontinuity and continuity}\label{ss.semi}

The quantity $\phi_p(S)$ can be used to give a short proof that $\pc(G)$ is lower semicontinuous in the local topology of transitive graphs: that is, $\liminf\pc(G_n)\geq \pc(G)$ holds; see \cite[Section 1.2]{DCT}. 
It can be proven for transitive graphs as follows: let $p<p_c(G)$, let $S\subset G$ be a set with $\phi_p^G(S)<1$ and let $r$ be such that $S\subset B_G(o,r)$. For $n$ large enough $B_{G_n}(o,r)\simeq B_G(o,r)$, hence $\phi_p^{G_n}(S)<1$, which implies $p\leq p_c(G_n)$. 
For bounded degree unimodular graphs, we will now show in a similar way that this inequality also holds for $\pe$; however, it fails for $\pl=\pc$, in general.

\begin{prop}\label{liminfpe}
Let $G_n$ and $G$ be unimodular random graphs with uniformly bounded degrees.
If $G_n$ converges to $G$ then $\liminf_{n\to\infty}\pe(G_n)\geq\pe(G)$.
\end{prop}

\begin{pf}
Let $p<\pe(G)$ and let $r$ be such that $\ee_{G}\left(\phi_p^\om(B_\om(o,r))\right)<1-\eps$ with some $\eps>0$.
Let $n$ be large enough to satisfy
\[ \sum_{H\in \mathcal{H}_{r+1}}
|\mu_{G_n}\left( B_\om(o,r+1)=H \right)-\mu_{G}\left( B_\om(o,r+1)=H \right)|<\frac{\eps}{2D^{r+1}},  \]
where $D$ is a uniform bound on the degrees of $G_n$ and $G$ and 
$\mathcal{H}_r$ is the set of possible $r$-neighbourhoods of the root in graphs with maximum degree $D$. 
Any $H\in \mathcal{H}_{r+1}$ satisfies $\phi_p^H(B_\om(o,r))\leq D^{r+1}$. We obtain 
\begin{align*}
\ee_{G_n}\left(\phi_p^\om(B_\om(o,r))\right)&
= \sum_{H\in \mathcal{H}_{r+1}}
\mu_{G_n}\left( B_\om(o,r+1)=H \right)\phi_p^H(B_\om(o,r)) \\
 &\leq \sum_{H\in \mathcal{H}_{r+1}}
\big\{ \mu_{G}\left( B_\om(o,r+1)=H \right)\phi_p^H(B_\om(o,r))\\
 &\quad + |\mu_{G_n}\left( B_\om(o,r+1)=H \right)-\mu_{G}\left( B_\om(o,r+1)=H \right)| \cdot |\ebd B_H(o,r)| \big\}\\
 &\leq \ee_{G}\left(\phi_p^\om(S)\right) + \frac{\eps}{2} <1. 
\end{align*}
It follows that $\pe(G_n)\geq p$ thus $\liminf \pe (G_n)\geq \pe(G)$. 
\end{pf}
\medskip

Now we prove Proposition \ref{liminf}, which states that if $G_n$ converges to a uniformly good unimodular graph $G$ in a uniformly sparse way, then $\pc(G_n)\to \pc(G)$. 
After the proof we present an example that shows how this proposition can be applied. 
Another application of the proposition appears in Example~\ref{product}. 

\begin{pfo} {\it Proposition \ref{liminf}.}
First, $G \subseteq G_n$ implies that $\pc(G)\geq \pc(G_n)$ for all $n$. 
For the sake of simplicity, we prove the inequality $\lim \pc(G_n)\geq \pc(G)$ for $k=1$. 
It can be proved for general $k$ in a similar way. 
Let $p<\pc(G)$. Our aim is to find a subset $B_n\in \fss(G_n)$ for $n$ large enough with $\phi_p^{G_n}(B_n)<1$. 
Let $n$ be sufficiently large to satisfy $r_n/2>R(p)$ and $c^{r_n/2}<\frac{1}{3}$. 
Fix a pair $(\om, \om_n)$ that satisfies the sparseness condition for $r_n$. 
Then, in the smaller ball $B_{\om_n}(o,r_n/2)$, there is at most one edge $\{x,y\} \in \om_n\setminus \om$. If this edge exists, let $B_n:=B_{\om_n}(o,r_n/2)\cup B_{\om_n}(x,r_n/2)\cup B_{\om_n}(y,r_n/2)$; otherwise, just let $B_n:=B_{\om_n}(o,r_n/2)$.  Note that $B_n\subset B_{\om_n}(o,r_n)$. Similarly, let $B:=B_{\om}(o,r_n/2)\cup B_{\om}(x,r_n/2)\cup B_{\om}(y,r_n/2)$, omitting those terms in the union that do not exists in $\omega$. (Note that it may happen that $x$ or $y$ does not exist in $\omega$, but not both, since $B_{\omega_n}(o,r_n/2)$ is connected.) 
The sets $B_n$ and $B$ satisfy $\ebd B_n=\ebd B$. 
We claim that we have $\phi_p^{\om_n}(B_n)<1$. 
There are three possibilities in terms of the edge $\{x,y\}$ for an open path connecting $o$ and a vertex $e^-$ in $B_n$: 
it connects $o$ and $e^-$ in $B$ or it connects $x$ or $y$ to $e^-$ in $B$. 
It follows that 
\begin{align*}
\phi_p^{\om_n}(B_n) &= p\sum_{e\in\ebd B_n}\pr^{\om_n}\left(o\xlra[B_n]{\om_n,p}e^-\right) \\
%& =p\sum_{e\in\ebd B_n} \left[ \pr^{\om_n}\left(o\xlra[B_n]{\om_n\setminus \{x,y\},p}e^-\right)
%+ \pr^{\om_n} \left( \{o\xlra[B_n]{\om_n,p}x\} \square \big\{\{x,y\}\textrm{ open}\big\} \square \{y\xlra[B_n]{\om_n,p}e^-\} \right) \right. \\
%& \hskip 1.3 in \left. +~ \pr^{\om_n} \left( \{o\xlra[B_n]{\om_n,p}y\} \square \big\{\{x,y\}\textrm{ open}\big\} \square \{x\xlra[B_n]{\om_n,p}e^-\} \right) \right]. 
%\end{align*}
%This can be bounded from above using the BK-inequality and the fact that the set $B_n$ can be changed to $B$ in each occurence in \eqref{Bxy} without altering the probabilies. We get  
%\begin{align*}
%& \leq p \sum_{e\in\ebd B} \left[ \pr^{\om}\left(o\xlra[B]{\om,p}e^-\right)
%+ p\,\pr^{\om} \left( \{o\xlra[B]{\om,p}x\} \square \{y\xlra[B]{\om,p}e^-\} \right)\right. \\
%& \hskip 1.75 in \left. +~ p\,\pr^{\om} \left( \{o\xlra[B]{\om,p}y\} \square \{x\xlra[B]{\om,p}e^-\} \right) \right] \\
& \leq p\sum_{e\in\ebd B} \left[ \pr^{\om}\left(o\xlra[B]{\om,p}e^-\right)
+ \pr^{\om}\left(y\xlra[B]{\om,p}e^-\right) + \pr^{\om}\left(x\xlra[B]{\om,p}e^-\right)\right] \\
& = \phi_p^{\om}(B) + \phi_p^{\om,y}(B) + \phi_p^{\om,x}(B) ~<~ 1
\end{align*}
by Lemma \ref{unifgoodexp}. 
If $x$ or $y$ does not exist in $\omega$, all its appearances in the above formulas involving $\omega$ can be replaced by the other vertex, and the inequalities remain true. 
It follows that $p\leq \pl(G_n)=\pc(G_n)$.
\end{pfo}

\begin{exa1}\label{pmexa}
The following example is a graph sequence $G_n$ where Proposition~\ref{liminf} applies. 
Let $G$ be a uniformly good unimodular graph of bounded degree; e.g., a unimodular quasi-transitive graph. % or a unimodular random tree of uniform sub-exponential growth. 
Let $H_n\subset V(G)$ be an invariant subset (i.e., given by a unimodular labelling) such that $\min\{\dist_G(x,y): x,y\in H_n\}\geq n$ almost surely. 
Such a subset can be produced as a factor of iid process: let $\{\xi_x: x\in V(G)\}$ be iid uniform random variables on $[0,1]$ and let 
$H_n:=\{x: \xi_x=\min\{\xi_y: y\in B_G(x,n)\}\}$. 
Consider now an invariant perfect matching of the points of $H_n$ (that is, an invariant partition of $H_n$ into pairs) and let $G_n$ be the union of that matching and $G$. An example of such
a perfect matching can be constructed as follows. Let $\{\zeta_e: e\in V(G)\}$ be iid uniform random variables on $[0,1]$ and 
consider the distance function $d$ on $V(G)$ defined as $d(x,y)=\inf \sum_{e\in P} \zeta_e$, where $P$ ranges over all paths connecting $x$ and $y$. It is easy to check that the infimum exists and is in fact a minimum; also, one can show that with the resulting metric the set $H_n$ is discrete, non-equidistant, and has no descending chains (see \cite{HPPS} for the definitions). By a method similar to the proof of Proposition 9 in \cite{HPPS}, one can show that the stable matching on $H_n$ is a perfect matching, just as desired.

For quasi-transitive graphs $G$, we have $p_T=p_c$. Then it is not surprising that, for any $p<p_c$, once $n$ is large enough, adding the sparse perfect matching cannot glue too many of the rather small finite clusters of $G$ together, and hence we still have $p<p_c(G_n)$. That is, one expects $p_c(G_n)\to p_c(G)$. This indeed holds by our general proposition, while an actual direct proof would need to handle some non-trivial technicalities.
\end{exa1}

Next, we turn to unimodular trees of uniform subexponential growth (see Definition~\ref{subexp.def}), proving Proposition~\ref{subexp}. This proposition gives further examples of uniformly good unimodular graphs (see Definition~\ref{unifgood}), while the convergence part will be used in Section~\ref{s.cost}.

\begin{pfo} {\it Proposition \ref{subexp}.} We start by proving the statement about the sequence $G_n$ with girth tending to infinity. By the uniform subexponential growth, for each $p<1$ there are positive integers $r=r(p)$ and $n_0(p)$ such that
\begin{align}\label{subexpvol}
|B_{G_n}(o_n,r)|\,p^r<1 
\end{align}
for every $n\geq n_0(p)$, almost surely. Now, by the girth tending to infinity, there exists $n_1(p)\geq n_0(p)$ such that, for every $n\geq n_1(p)$, the ball $B_{G_n}(o_n,r)$ is a tree, and therefore
\begin{align}\label{subexpphi}
\phi^{G_n}_p(B_{G_n}(o_n,r)) \leq |B_{G_n}(o_n,r)|\,p^r\,. 
\end{align}
Combining (\ref{subexpvol}) and  (\ref{subexpphi}), and taking $p\to 1$, the balls $B_{G_n}(o_n,r)$ show that $\pl(G_n)$ and $\pe(G_n)$ tend to 1. By Theorem~\ref{pl=pc}, we also have $\pc(G_n)\to 1$.

Now, if $G$ is a unimodular tree of subexponential growth, then~(\ref{subexpphi}) holds for every $r$, hence $\pe(G)=\pl(G)=\pc(G)=1$, and uniform goodness is also clear from the definition. Then Corollary~\ref{pc=pt} implies $\pt(G)=\ptk(G)=1$, as well.
\end{pfo}

\subsection{Counterexamples}\label{ss.loc.countex}

Our first example will show that even if we keep the condition of uniformly sparse convergence of $G_n$ to $G$ of Proposition~\ref{liminf}, without $G$ being uniformly good, the conclusion may not hold. Next, Example~\ref{Zd} will show that keeping the limit uniformly good but removing the condition of uniform sparseness will make the conclusion false. Finally, Example~\ref{box} will show that the inequality of the lower semicontinuity may be strict even when invariant subgraphs $G_n$ of $\Z^2$ converge to $\Z^2$.

\begin{exa}\label{lamplighter}
There exists a sequence $(G_n)$ of invariant random subgraphs of a Cayley graph, converging to an invariant subgraph $G$ in a uniformly sparse way, such that $\lim\pc(G_n)< \pc(\lim G_n)$.
\end{exa}

\begin{pf} The first step is to construct an invariant percolation on a Cayley graph of the lamplighter group all whose clusters are isomorphic to the canopy tree $\Lambda$ (see Definition \ref{canopytree}. In more detail:

Consider the generators  $\{Rs, R, sL, L\}$ of the lampligher group $\Z_2 \wr \Z = \oplus_{\Z} \Z_2 \rtimes \Z$, where $R:=(0,1), L:=(0,-1)$, and $s:=(e_0,0)\in \Z_2 \wr \Z$ with $e_0\in \{0,1\}^\Z$, $(e_0)_j=\delta_{0,j}$. It is well-known (see, e.g., \cite{W}) that the Cayley graph with respect to these generators is the Diestel-Leader graph DL(2,2). This graph can be defined using two trees $\T^1$ and $\T^2$ which both are 3-regular infinite rooted trees with a distinguished end and {Busemann functions} $\mathfrak{h}_i:\T^i\to \Z, i=1,2$, as in Definition \ref{canopytree}. Each vertex $x\in\T^i$ has exactly one neighbour $\bar{x}$ with $\mathfrak{h}_i(\bar{x})=\mathfrak{h}_i(x)-1$, called the parent of $x$. We call the other two neighbours the children of $x$. Now consider the following percolation on $\T^1$: for each vertex $x$ we delete the edge connecting $x$ to one of its two children, independently with equal probabilities.
We get a random subgraph of $\T^1$ consisting of infinite simple paths.  We then delete the edges in the graph DL(2,2) whose first coordinate is a deleted edge in $\T^1$. The resulting random subgraph $\F\subset$ DL(2,2) is invariant under the action of the lamplighter group and it consists of infinitely many components which are all isomorphic to the canopy tree $\Lambda\subset\T$. The probability that the root is in the $n^\textrm{th}$ level of its component in $\F$ is clearly $2^{-n-1}$. The canopy tree with a random root chosen according to this distribution is a unimodular random graph, as it also must be the case by Proposition~\ref{caysub}.

The significance of the canopy tree for this construction (as in Example~\ref{subexpcanopy}) will be that it has one end, thus $\pc(\Lambda)=1$, while one can easily compute that $\pt(\Lambda)={1}/{\sqrt{2}}$.

Now let $\GG$ be the free product of $\Z_2:=\Z/{2\Z}$ and the lamplighter group $\Z_2 \wr \Z$.
Let $\Gamma$ be the left Cayley graph of $\GG$ with respect to the generators $\{a, Rs, R, sL, L\}$ 
where $a$ is the generator of the free factor $\Z_2$. Let  $\beta:\GG\longrightarrow \Z_2\wr\Z$ be the natural projection homomorphism:
if $w=a_1b_1\dots a_kb_k$ is a word in $\GG$ such that $a_j\in\Z_2, b_j\in\Z_2\wr\Z, j=1,\dots,k$, then $\beta(w):=b_1\dots b_k\in \Z_2\wr\Z$. We now define $G$ to be the following random spanning subgraph of $\Gamma$: 
let $e$ be in $E(G)$ if{f} $\beta(e^-)$ and $\beta(e^+)$ are connected by an edge in $\F$. 
The distribution of $G$ is invariant under the action of $\GG$ and each component of $G$ is a canopy tree, 
hence $\pc(G)=1$.

We define a sequence $(G_n)$ of random subgraphs of $\Gamma$ converging to $G$. 
We choose an element $b\in\{0,1,\dots n-1\}$ uniformly at random. 
For each vertex in $L_{\T^1}(b+kn), k\in\Z$ we choose 
one of its descendants in $L_{\T^1}\left(b+(k+1)n\right)$ uniformly at random and 
we choose all vertices in $L_{\T^2}\left(-b+kn\right)$. 
Let $S_n$ be the set of edges $e\in E(\Gamma)$ such that $e$ is labelled by the generator $a$ and 
both coordinates of $\beta(e^-)=\beta(e^+)$ are chosen vertices in the above procedure.  
Let $G_n:=G\cup S_n$.

We show that $\pc(G_n)\leq \frac{1}{\sqrt{2}}$ for all $n$. 
Let $p>\frac{1}{\sqrt{2}}=\pt(\Lambda)$, let $n$ be a positive integer and 
consider Bernoulli($p$) percolation on $G_n$. 
Denote by $T(v)$ the component of the vertex $v$ in $G$ and 
by $\C_v$ the component of the vertex $v$ in the percolation on $G_n$. 
Let $s(v):=\min\{l: L_{T(v)}(l)\cap S_n\neq\emptyset\}$. 
We define a branching process depending on the percolation on $G_n$. 
For each vertex $v$ of $\Gamma$ let 
$N_v:=\{ax: x\in T(v)\cap \C_v\cap S_n\setminus\{v\}, \{x,ax\}\textrm{ is open}\}$. 
Let $Z_1:=N_o$ and let $Z_{k+1}:=\bigcup_{v\in Z_k} N_v$.  
Note that $Z_i\neq Z_j, i\neq j$ and $Z_j\subset \C_o$. 
The distribution of $|N_v|$ depends only on the level of $v$ in $T(v)$ and on $s(v)$. 
The distribution of $|N_v|$ conditioned on $\{o\in L_{T(o)}(l), s(v)=s\}$ with any $l$ and $s$ stochastically dominates 
the distribution of $|N_v|$ conditioned on the event $\{v\in L_{T(v)}(0), s(v)=n-1\}$. 
Therefore the distribution of $|Z_k|$ stochastically dominates the distribution of 
the $k^\textrm{th}$ generation of the Galton--Watson process 
with offspring distribution $|N_v|$ conditioned on $\{v\in L_{T(v)}(0), s(v)=n-1\}$, which has infinite expectation. 
Hence $\pp(\liminf |Z_k|>0)>0$ 
which implies $\pp(|\C_o|=\infty)>0$.  
\end{pf}

\begin{exa}\label{Zd}
There exists a sequence $(G_n)$ of invariant random subgraphs of a Cayley graph such that 
$\lim\pc(G_n)< \pc(\lim G_n)$ and $\lim G_n$ is uniformly good.
\end{exa}

\begin{pf} Let $\Gamma$ be a Cayley graph of a finitely generated group such that there exists a random subgraph $\bar{G}$
which satisfies the following: the distribution of $\bar{G}$ is invariant under the action of the group,
it consists of infinitely many infinite components and
each component has critical percolation probability $\bar{p}<1$. (A very simple example is that $\Gamma$ is $\Z^d$ and $\bar{G}$ is a lamination by copies of $\Z^{d-1}$, with $d\geq 3$.) Let $G'$ be an invariant random connected subgraph of $\Gamma$ such that
$\pc(G')>\bar{p}$. For example, if $\Gamma$ is amenable, then one can choose $G'$ to be an
invariant spanning tree of $\Gamma$, which always exists and has at most two
ends, and hence $\pc(G')=1$; see \cite{BLPS}, Theorem 5.3. Moreover, if $\Gamma$ has sub-exponential
volume growth (see Definition \ref{subexp.def}), then so does the spanning tree $G'$, and it is uniformly good by Proposition~\ref{subexp}.

Now let $\eps_n\to 0$ be a sequence of positive numbers and let $G_n$ be the following random subgraph of $\Gamma$:
we remove each component of $\bar{G}$ with probability $1-\eps_n$ and keep it with probability $\eps_n$
independently for each component.
Let $G_n$ be the union of $G'$ and the remaining components of $\bar{G}$. 
It follows from Proposition~\ref{caysub} that $G_n$ is unimodular.
The sequence $(G_n)$ converges to $G'$, but $\pc(G_n)\leq\bar{p}<\pc(G')$ for each $n$. The sequence $p_c(G_n)$ has a convergent subsequence, hence we can choose the corresponding subsequence $\eps_{k(n)}$, and get $\lim \pc(G_{k(n)})\leq \bar{p}<\pc(G')$. 

We get a similar example that is uniformly good if we set $\Gamma:=\Z^5$,
$\bar{G}:=\bigcup_{y\in\Z^2}\{y\}\times \Z^3$ and
$G':=\bigcup_{x\in\Z^3}\Z^2\times\{x\}$. In this example $G'$ is not
connected, but each $G_n$ is connected almost surely, and $\pc(G_n)\leq
\pc(\Z^3)<\pc(\lim G_n)=\pc(\Z^2)<1$ for each $n$. 
\end{pf}

\begin{exa}\label{box}
There exists a sequence $(G_n)$ of invariant random subgraphs of a Cayley graph such that 
$1>\lim\pc(G_n)> \pc(\lim G_n)$.
\end{exa}

\begin{pf} We define $G_n$ as a vertex and edge replacement (see Subsection~\ref{ss.opera} and \cite{AL07}, Example 9.8) 
of $\Z^2$ where we replace each vertex $x$ by the graph $Q_x$ isomorphic to 
$Q_n$ and we replace each edge by a path of length two that 
joins the middle points of the neighbouring sides of the boxes corresponding to the endpoints of the edge. 
The graphs $G_n$ can be considered as deterministic subgraphs of $\Z^2$ with a randomly chosen root. 
The sequence $G_n$ converges to $\Z^2$. 

We show that $\frac{1}{2}<\lim\pc(G_n)<1$. 
Denote by $G_n(p)$ the subgraph obtained by the Bernoulli($p$) percolation on $G_n$, and let $H_n(p)$ be the following percolation on $\Z^2$: let an edge $\{x,y\}$ open, iff both edges are open in the path that joins the boxes $Q_x$ and $Q_y$ in $G_n$. 
The existence of an infinite cluster in $G_n(p)$ implies the existence of an infinite cluster in $H_n(p)$. The law of $H_n$ equals the law of the Bernoulli($p^2$) percolation on $\Z^2$, hence $\pc(G_n)\geq \frac{1}{\sqrt{2}}$ for each $n$. 

To show that $\limsup p_c(G_n)<1$, we define the percolation $\bar{H}_n(p)$ on $\Z^2$. Denote by $\mathcal{A}_x(n)$ the event that the vertices $(0,-n)$, $(0,n)$, $(-n,0)$, $(n,0)$ are in the same cluster in Bernoulli($p$) percolation on the box $Q_x\subset G_n$. Let an edge $\{x,y\}\in \bar{H}_n(p)$, iff $\{x,y\}\in H_n(p)$, and both of the events $\mathcal{A}_x(n)$ and $\mathcal{A}_y(n)$ occurs. The existence of an infinite cluster in $\bar{H}_n(p)$ implies the existence of an infinite cluster in $G_n(p)$.
Let $1>p_0>\frac{1}{2}$ be arbitrary. There is an $\eps>0$ such that if the marginals of a 2-dependent percolation on $\Z^2$ are at least $(1-\eps)^4$, then this percolation stochastically dominates Bernoulli($p_0$) percolation; see \cite[Theorem 0.0]{LSS}. 
Lemma \ref{crossing} implies, that we can find constants $1-\eps<p_1<1$ and $N$ such that for any $p>p_1$, $n\geq N$ and for any vertex $x\in V(\Z^2)$ the event $\mathcal{A}_{x}(n)$ occurs whith probability at least $1-\eps$, thus $\pr(e\in \bar{H}_n(p))\geq p_1^2(1-\eps)^2\geq (1-\eps)^4$ for any edge $e\in E(\Z^2)$. The events $\{e_1\in\bar{H}_n\}$ and $\{e_2\in\bar{H}_n\}$ are independent if the distance of $e_1$ and $e_2$ is at least 2, hence $\bar{H}_n(p)$ stochastically dominates Bernoulli($p_0$) percolation. It follows that $\limsup p_c(G_n)\leq p_1<1$.
\end{pf}

\section{On transitive graphs of cost 1}\label{s.cost}

As proved in \cite[Theorem 5.3]{BLPS}, a transitive graph $G$ is amenable if and only if it has an invariant spanning tree $\TT$ with at most two ends, hence with expected degree $2$ and $\pc(\TT)=1$. Briefly: for the existence of $\TT$ for an amenable $G$, see the proof of Corollary~\ref{biLip} below, while from an invariant connected spanning graph $\TT$ with $\pc(\TT)=1$ it is not hard to construct an invariant mean on $G$, and thus deduce amenability.

Proposition~\ref{subexp} tells us that, under the stronger condition of subexponential growth, we get a spanning tree $\TT$ with the stronger property $\pt(\TT)=\ptk(\TT)=1$. Moreover, we can achieve approximately 1-dimensional percolation behaviour $\pc(G_k)\to 1$ via connected spanning subgraphs that have the same large-scale geometry as $G$. 

\begin{cor}\label{biLip}
If $G$ is a transitive amenable graph, then there is a
sequence of invariant random subgraphs $G_k$ which satisfies the following:
each $G_k$ is a bi-Lipschitz (in particular, connected) spanning subgraph of $G$,
the girth of $G_k$ tends to infinity and 
$G_k$ locally converges to an invariant random spanning tree $\TT$ with at most two ends.

If $G$ is a transitive graph with sub-exponential volume growth, then $\lim \pc(G_k)=1$.
\end{cor}

\begin{pf}
We construct $\TT$ as in \cite{BLPS}, Theorem 5.3:
let $F_n$ be a sequence of F\o lner sets such that $\sum_{n=1}^{\infty}\frac{|\ebd F_n|}{|F_n|} <1$.
For each $n$ and $x\in V(G)$ choose a random $g_{x,n}\in Aut(G)$ that takes $o$ to $x$, and
a random bit $Z_{x,n}$ that equals 1 with probability $\frac{1}{|F_n|}$. Choose all $g_{x,n}$ and $Z_{x,n}$ independently.
Let $\om_n:=E(G)\setminus \bigcup_{x\in V(G), Z_{x,n}=1}\ebd(g_{x,n}F_n)$; i.e.,
we remove all edges in the boundaries of the translates of $F_n$ with $Z_{x,n}=1$.
Let $\bar{\om}_n=\bigcap_{k\geq n} \om_k$. Each $\bar{\om}_n$ has only finite components.

To construct $\TT$ and $G_k$, choose uniform labels $L_e$ in [0,1] independently for each $e\in E(G)$.
For each finite component of $\bar{\om}_1$ take the minimal spanning tree of the component with respect to the labels.
Denote by $T_1$ the union of these trees.
Let $T_2$ be the union of $T_1$ and the edges in $\bar\om_2\setminus\bar\om_1$ with minimal labels such that
the components of $T_2$ are spanning trees of the components of $\bar\om_2$.
Continue inductively, and let $\TT:=\bigcup T_n$. This is an invariant random spanning tree, which has at most 2 ends (otherwise it would have infinitely many ends, which is impossible, since $G$ is amenable).

To construct $G_k$ we define a color for each edge. Let all edges in $\TT$ be green.
In each component of $\bar\om_1$ do the following:
consider the edge with the smallest label which has no color.
If there is a path of length at most $k$ between its endpoints consisting of green edges, then color it red, otherwise color it green.
Continue inductively for the edges in the component. This procedure defines a color for each edge of $\bar\om_1$.
If all edges in $\bar\om_n$ have a color, then
continue coloring the edges of $\bar\om_{n+1}\setminus\bar\om_n$ in the same way.
Let $G_k$ be the union of the green edges.
It follows from the construction that $G_k$ is invariant, its girth is at least $k+2$ and
for each edge of $G$ there is a path in $G_k$ between its endpoints with length at most $k$.
The sequence $G_k$ converges to $\TT$.

If $G$ has sub-exponential volume growth, then so does $\TT$ and each $G_k$, and all of them are unimodular (by~\cite[Corollary 1]{SW} and Proposition~\ref{caysub} above). Thus $p_c(G_k)\to 1$ follows from Proposition~\ref{subexp}.
\end{pf}
\medskip

It might be surprising at first sight that, as opposed to having a spanning subgraph with $\pc=1$, the existence of a sequence $G_k$ as in the corollary does not imply amenability:

\begin{exa}\label{product}
$\T_3 \times \Z$ has a sequence of invariant bi-Lipschitz subgraphs $G_k$ with $\pc(G_k)\to 1$. 
\end{exa}

\begin{pf}
One can partition the edges of $\T_3$ into 3 disjoint perfect matchings $M_1$, $M_2$ and $M_3$ in an invariant way. (See, for instance, \cite{LyFIID}, around Proposition 2.4.) Then, consider the following subgraphs $G_k\subseteq \T_3\times \Z$: we keep all the edges in the subgraphs $\{v\}\times\Z$ and the edges $\{e\}\times \{3jk+ik\}$ where $e\in M_i$, $j\in \Z$. We choose a uniform random integer $b\in\{0,\dots k-1\}$ and translate this subgraph by $(id,b)$ to get the invariant subgraph $G_k$ of $\T_3\times\Z$. 
Each $G_k$ is clearly bi-Lipschitz equivalent to $\T_3\times \Z$. On the other
hand, we have $\pc(G_k)\to 1$: either from Proposition~\ref{liminf}, or more directly, by observing that the universal cover $T_k$ of $G_k$ can be obtained from $\T_3$ by replacing ``two thirds'' of the edges by a path of length $k$; for this tree, it is easy to see that $\pc(T_k)\to 1$, while $\pc(T_k)\leq\pc(G_k)$ holds by \cite[Theorem 6.47]{LP}. 
\end{pf}
\medskip

So, what is the class of transitive graphs for which the existence of such a sequence $G_k$ may be expected? The answer seems to have something to do with the notion of cost from measurable group theory. (See Subsection~\ref{ss.results} for references.) The {\it cost of a group} $\GG$ is defined as half of the infimum of the expected degrees of its invariant connected spanning graphs. The {\it $\GG$-cost of a transitive graph} $G$ may be defined similarly, over $\GG$-invariant random connected spanning {\it sub}graphs of $G$, where $\GG\leq Aut(G)$ is a vertex-transitive subgroup of graph-automorphisms. It is not known in general that, if we first fix a Cayley graph $G$ of $\GG$, then the $\GG$-cost of $G$ is always as small as the cost of $\GG$ (which is the cost of the complete graph on $\GG$). Nevertheless, we have seen that cost 1 can be achieved inside any Cayley graph of any amenable group (since the expected degree of an infinite unimodular tree with at most two ends is 2). 

We will now show that a sequence of invariant spanning subgraphs $G_k$ with $\pc(G_k)\to 1$ implies that the cost is 1. The bi-Lipschitz condition does not appear here, but it is quite possible that once we have a sequence with $\pc(G_k)\to 1$, it can always be modified to fulfill the bi-Lipschitz property, as well. Note that the bi-Lipschitz condition is also natural from the point of view of Elek's combinatorial cost for sequences of finite graphs \cite{E}.

\begin{lemma}\label{cost}
If $\Gamma$ is a Cayley graph of $\GG$, and there exists a sequence of $\GG$-invariant connected spanning subgraphs $G_k \subset \Gamma$ with $\pc(G_k)\to 1$, then the cost of $\Gamma$, hence of $\GG$, is 1.
\end{lemma}

\begin{pf} 
Take $\eps_k\to 0$ such that $\pc(G_k)>1-\eps_k$. Then, all clusters of Bernoulli($1-\eps_k$) percolation on $G_k$ are finite almost surely. Let the set of closed edges be denoted by $\eta_k \subset G_k \subset \Gamma$, an invariant percolation itself. In each finite cluster, take a uniform random spanning tree, a subtree of $G_k$. The union of all these finite spanning trees and $\eta_k$ will be $\om_k$. One the one hand, it is clear that $\om_k$ is a connected spanning subgraph of $G_k$, hence of $\Gamma$. On the other hand, the expected degree of $o$ in $\om_k$ is at most $\ex \deg_{\eta_k}(o) + 2 \leq d\eps_k+ 2$, where $\deg_\Gamma(o)=d$. As $k\to\infty$, we obtain that the cost of $\Gamma$ is 1.
\end{pf}
\medskip

We do not know if the converse of Lemma~\ref{cost} holds: 

\begin{quest}
Does there exist,  for any Cayley graph $G$ of any group $\GG$ of cost 1, a sequence of $\GG$-invariant bi-Lipschitz spanning subgraphs $G_k \subset G$ with $\pc(G_k)\to 1$? At least for amenable $G$?
\end{quest}

For amenable Cayley graphs $G$, a first step of independent interest could be a positive answer to the following question, mentioned in Subsection~\ref{ss.results}:

\begin{quest}
For any amenable Cayley graph, is there an invariant random spanning subtree of subexponential growth? More boldly, does there always exist an invariant random Hamiltonian path? 
\end{quest}

\section*{Acknowledgments}

We are grateful to S\'ebastien Martineau and two referees for comments on the manuscript.

Our work was partially supported by the ERC Consolidator Grant 648017
(DB), the Hungarian National Research, Development and Innovation
Office, NKFIH grant K109684 (GP and \'AT), an MTA R\'enyi Institute
``Lend\"ulet'' Research Group (GP), and an EU Marie Curie Fellowship
(\'AT).

\ \\
\noindent
{\bf Dorottya Beringer}\\
Alfr\'ed R\'enyi Institute of Mathematics, Hungarian Academy of Sciences\\
Re\'altanoda u. 13-15, Budapest 1053 Hungary\\
and Bolyai Institute, University of Szeged\\
Aradi v\'ertan\'uk tere 1, Szeged 6720 Hungary\\
\texttt{beringer[at]math.u-szeged.hu}\\
\ \\
{\bf G\'abor Pete}\\
Alfr\'ed R\'enyi Institute of Mathematics, Hungarian Academy of Sciences\\
Re\'altanoda u. 13-15, Budapest 1053 Hungary\\
and Institute of Mathematics, Budapest University of Technology and Economics\\
Egry J. u. 1., Budapest 1111 Hungary\\
\url{http://www.math.bme.hu/~gabor}\\
\ \\
{\bf \'Ad\'am Tim\'ar}\\
Alfr\'ed R\'enyi Institute of Mathematics, Hungarian Academy of Sciences\\
Re\'altanoda u. 13-15, Budapest 1053 Hungary\\
\texttt{madaramit[at]gmail.com}\\

\end{document}